\documentclass[preprint]{revtex4}

\usepackage[ansinew]{inputenc}
\usepackage{amsmath, amsfonts} 
\usepackage{amsthm, amssymb}
\usepackage[english]{babel}
\usepackage[T1]{fontenc}
\usepackage{lmodern}
\usepackage{graphicx}
\usepackage{subfigure}

\theoremstyle{definition}

\newtheorem{remark}{Remark}

\usepackage{verbatim}
\usepackage{array}
\usepackage{float}
\usepackage{subfigure}
\usepackage{color}

\newcommand{\guillemets}[1]{``#1''}
\newcommand{\ve}[1]{\mathbf{#1}}

\begin{document}

\title{Global computation of phase-amplitude reduction for limit-cycle dynamics}
\author{A. Mauroy}
\email{alexandre.mauroy@unamur.be}
\affiliation{Department of Mathematics and Namur Institute for Complex Systems, University of Namur, B-5000 Namur, Belgium}
\author{I. Mezic}
\email{mezic@engr.ucsb.edu}
\affiliation{Department of Mechanical Engineering, University of California Santa Barbara, Santa Barbara, CA 93106, USA}

\begin{abstract}
Recent years have witnessed increasing interest to phase-amplitude reduction of limit-cycle dynamics. Adding an amplitude coordinate to the phase coordinate allows to take into account the dynamics transversal to the limit cycle and thereby overcomes the main limitations of classic phase reduction (strong convergence to the limit cycle and weak inputs). While previous studies mostly focus on local quantities such as infinitesimal responses, a major and limiting challenge of phase-amplitude reduction is to compute amplitude coordinates globally, in the basin of attraction of the limit cycle.

In this paper, we propose a method to compute the full set of phase-amplitude coordinates in the large. Our method is based on the so-called Koopman (composition) operator and aims at computing the eigenfunctions of the operator through Laplace averages (in combination with the harmonic balance method). This yields a forward integration method that is not limited to two-dimensional systems. We illustrate the method by computing the so-called isostables of limit cycles in two, three, and four-dimensional state spaces, as well as their responses to strong external inputs.

\end{abstract}

\maketitle

\noindent\textbf{Oscillatory behaviors in biology, physics, and engineering are often related to high-dimensional limit-cycles dynamics. These possibly complex, high-dimensional dynamics can be reduced to simple, low-dimensional dynamics through phase reduction. While classic phase reduction only captures the effect of small perturbations to the system, recent developments have introduced a more general phase-amplitude reduction, which is well-suited to large perturbations. This reduction is related to phase-amplitude coordinates associated with specific families of sets in the state space: the isochrons and the isostables. As a main limitation of phase-amplitude reduction, the computation of isostables is intricate and typically limited to local quantities such as infinitesimal responses. This paper presents a numerical method to compute the isostables in the whole basin of attraction of the limit cycle. This method relies on the framework of the Koopman operator, which allows to interpret the isostables as level sets of specific eigenfunctions.}

\section{Introduction}

High-dimensional limit cycle dynamics can be reduced to simple one-dimensional dynamics through phase reduction \cite{Malkin,Winfree,Kuramoto_book}. In this framework, the phase variable evolving on the circle describes the state of the system and the phase response indicates the effect of external inputs. This powerful reduction has proved to be an accurate and efficient tool to describe the response of complex oscillatory dynamics, in particular in the context of neuroscience, where it is also convenient from an experimental point of view and useful to study collective behaviors (see e.g. \cite{Brown, Ermentrout_book, Izi_book}).

However, phase reduction is not valid when the convergence rate toward the limit cycle is too slow or when external inputs are too strong. This is due to the fact that the phase response does not capture the full system dynamics, but only the dynamics in the neighborhood of the limit cycle. For this reason, the past years have witnessed increasing effort to overcome this limitation. Higher order approximations of phase responses were proposed in \cite{Suvak,Takeshita} and specific phase responses taking into account the effect of a train of several pulses were considered in \cite{Klinshov, Oprisan}, among others. Alternatively, an elegant approach consists in augmenting the phase space with an amplitude coordinate which takes into account the dynamics transversal to the limit cycle \cite{Llave_param,Guillamon, Nakao2013, Wedgwood}. In this case, one obtains a simple, reduced action-angle representation of the limit cycle dynamics. We focus on this approach in the present paper.

Phase-amplitude reduction is strongly connected to spectral properties of the so-called Koopman operator \cite{Budisic_Koopman, Mezic}. It was shown in \cite{Mauroy_Mezic} that a specific eigenfunction of the Koopman operator can be used to define the phase coordinate, or equivalently that the level sets of this eigenfunction are the isochrons of the limit cycle \cite{Winfree2}. In the case of systems with a stable equilibrium, the spectral properties of the Koopman operator were used to define an amplitude coordinate through a family of sets called isostables, which complements the family of isochrons \cite{MMM_isostables}. The approach based on isochrons and isostables was extended in a straightforward way to limit cycles \cite{Wilson_isostable} and used recently in the context of optimal control \cite{Monga}. This extension is in fact equivalent to phase-amplitude reduction proposed in \cite{Llave_param,Guillamon}, since both induce a constant rate of convergence toward the limit cycle in the reduced coordinates.

Since the goal of phase-amplitude reduction is to take into account large perturbations that drive the state away from the limit cycle, it is natural to compute amplitude coordinates, or equivalently isostables, in the large. This not only allows to fully characterize the sensitivity to large perturbations, but also provides a global picture of the limit cycle dynamics. However, the global computation of amplitude coordinates is delicate and, to the authors knowledge, previous contributions mainly focused on local quantities such as the infinitesimal isostable (or phase) response (see \cite{Guillamon2, Nakao_phase_amplitude} in the large and \cite{Wilson_isostable} along the limit cycle).

In this paper, we go a step further by computing the full set of phase-amplitude coordinates in the large. To do so, we exploit the Koopman operator framework and compute so-called Fourier and Laplace averages yielding the eigenfunctions of the operator. This method can be seen as an extension of the results of \cite{MMM_isostables} to limit cycles, although numerical computations are more involved in this case and based on a Fourier expansion of the limit cycle (e.g. through the harmonic balance technique \cite{harmonic_balance}). Our numerical method relies on forward integration and is efficient in high-dimensional systems. As shown in this paper, it can be used to compute the global isostables of a (four-dimensional) limit cycle. In addition to phase-amplitude coordinates, this forward-integration method can provide the (infinitesimal) phase and isostable responses, thereby complementing previous approaches based on the adjoint method and backward integration \cite{Guillamon2, Nakao_phase_amplitude}.

The rest of the paper is organized as follows. In Section \ref{sec:Koopman}, we introduce phase-amplitude reduction within the framework of the Koopman operator, showing that the reduced coordinates are directly related to the eigenfunctions of the Koopman operator. In Section \ref{sec:comput_eigen}, a method to compute the eigenfunctions of the Koopman operator (and their gradient) is presented, based on Fourier and Laplace averages estimated along the trajectories of the system. Section \ref{sec:num} provides a few guidelines for numerical computation and Section \ref{sec:examples} illustrates the method with examples in two, three, and four dimensional state spaces. Phase-amplitude coordinates are also used to study the effect of an external input on the system. Finally, concluding remarks are given in Section \ref{sec:conclu}.

\section{From Koopman operator to phase-amplitude reduction}
\label{sec:Koopman}

\subsection{Koopman operator}

Consider a dynamical system
\begin{equation}
\label{eq:syst}
\dot{\ve{x}}  = \ve{F}(\ve{x}) \qquad \ve{x} \in \mathbb{R}^n
\end{equation}
which is equivalently described by the flow $\varphi: \mathbb{R}^+ \times \mathbb{R}^n \to \mathbb{R}^n$, i.e. $\varphi(t,\ve{x}_0) = \varphi^t(\ve{x}_0)$ is a solution to \eqref{eq:syst} with the initial condition $\ve{x}_0$. We assume that $\ve{F}$ is Lipschitz, so that this solution exists and is unique.

The group of Koopman operators associated with \eqref{eq:syst} is given by
\begin{equation*}
U^t:\mathcal{F} \to \mathcal{F}, \quad U^t f = f \circ \varphi^t,
\end{equation*}
for all $t \geq 0$ and all functions $f\in \mathcal{F}$, where $\mathcal{F}$ is a well-defined linear vector space that contains constant functions. A function $\phi_{\lambda} \in \mathcal{F}$ is an eigenfunction of the Koopman operator if it satisfies
\begin{equation*}
U^t \phi_\lambda = e^{\lambda t} \phi_{\lambda} \quad \forall t \geq 0
\end{equation*}
for some $\lambda\in \mathbb{C}$. The value $\lambda$ is the corresponding eigenvalue and belongs to the point spectrum of the operator. It is easy to see that the constant function is an eigenfunction of the operator associated with the eigenvalue $\lambda=0$. Other eigenfunctions and eigenvalues capture the dynamics of the underlying system. 

Assume now that the system \eqref{eq:syst} admits a limit cycle $\Gamma$ (of frequency $\omega$), with a basin of attraction $\mathcal{B}(\Gamma) \subseteq \mathbb{R}^n$. If the limit cycle is stable and normally  hyperbolic, i.e. its Floquet exponents $\Lambda_j$ satisfy
\begin{equation*}
\Lambda_0 = 0,\quad  \Re\{\Lambda_j\} < 0 \quad \forall j=1,\dots,n-1,
\end{equation*}
then the spectrum of the Koopman operator is completely characterized \cite{Mezic2017}. It includes the Floquet exponents and also captures the limit cycle frequency: the principal eigenvalues of the Koopman operator \cite{Mohr} are given by
\begin{equation*}
\lambda_1 = i \omega, \quad \lambda_j = \Lambda_{j-1} \quad j=2,\dots,n
\end{equation*}
and there exist associated eigenfunctions 
\begin{equation}
\label{eq:eigen_cycle}
\phi_{\lambda_1} = \phi_{i \omega}, \quad \phi_{\lambda_j} = \phi_{\Lambda_{j-1}} \quad j=2,\dots,n.
\end{equation}
that have support on $\mathcal{B}(\Gamma)$ and are continuously differentiable in the interior of $\mathcal{B}(\Gamma)$ \cite{Mauroy_Mezic_stability}.

\subsection{Phase-amplitude reduction of limit-cycle dynamics}

The eigenfunctions of the Koopman operator can be used to derive an appropriate set of linearizing coordinates. Consider the eigenfunction $\phi_{\lambda}$ and the new variable $z(t) = \phi_{\lambda}(\varphi^t(\ve{x}))$. Then , we have
\begin{equation*}
\dot{z}(t) = \frac{d}{dt} \left[\phi_{\lambda}(\varphi^t(\ve{x})) \right] = \frac{d}{dt} \left[U^t \phi_{\lambda}(\ve{x}) \right] = \frac{d}{dt} \left[ e^{\lambda t} \phi_{\lambda}(\ve{x}) \right] = \lambda \left[ e^{\lambda t} \phi_{\lambda}(\ve{x}) \right] = \lambda z(t).
\end{equation*}
Moreover, if the operator admits $n$ eigenfunctions $\phi_{\lambda_j}$, $j=1,\dots,n$, such that
\begin{equation}
\label{eq:h}
\ve{h}:X \to \mathbb{C}^n, \quad \ve{x} \mapsto \ve{h}(\ve{x})=(\phi_{\lambda_1}(\ve{x}),\dots,\phi_{\lambda_n}(\ve{x}))
\end{equation}
is a diffeomorphism on a set $X$, then the dynamics \eqref{eq:syst} on $X$ are given by
\begin{equation}
\label{eq:lin_dyn}
\dot{\ve{z}} = \begin{pmatrix} \lambda_1 & 0 & \cdots & 0 \\
0 & \lambda_2 & \cdots & 0 \\
\vdots & \vdots & \ddots & \vdots \\
0 & 0 & \cdots & \lambda_n \end{pmatrix} \ve{z}
\end{equation}
in the new coordinates $\ve{z}=\ve{h}(\ve{x})$. In the case of normally hyperbolic limit cycles, it is shown in \cite{Lan} that the transformation of coordinates \eqref{eq:h} with the eigenfunctions \eqref{eq:eigen_cycle} is a diffeomorphism from $X=\mathcal{B}(\Gamma)$ to $\mathbb{C}^n$. It follows that the limit-cycle dynamics can be described by \eqref{eq:lin_dyn}. Moreover, if $\lambda_j \in \mathbb{R}$, we define $r_j = z_j = \phi_{\lambda_j}$ and \eqref{eq:lin_dyn} yields $\dot{r}_j = \sigma_{j} \, r_j$ with $\sigma_j \triangleq \Re\{\lambda_j\} = \Re\{\Lambda_{j-1}\}$. If $\lambda_j \notin \mathbb{R}$, we define $r_j = |z_j| = |\phi_{\lambda_j}|$, $\theta_j = \angle z_j = \angle \phi_{\lambda_j}$ and \eqref{eq:lin_dyn} yields $\dot{r}_j = \sigma_{j} \, r_j$, $\dot{\theta}_j = \omega_j$ with $\sigma_j \triangleq \Re\{\lambda_j\} = \Re\{\Lambda_{j-1}\}$ and $\omega_j \triangleq \Im\{\lambda_j\} = \Im\{\Lambda_{j-1}\}$. Eliminating redundant variables (due to complex conjugate eigenvalues and eigenfunctions) and reordering the indices, we obtain
\begin{equation*}
\begin{array}{lcl}
\dot{\theta}_1 & = & \omega \\
\dot{\theta}_2 & = & \omega_2 \\
& \vdots & \\
\dot{\theta}_{m+1} & = & \omega_{m+1} \\
\dot{r}_2 & = & \sigma_{2} \, r_2 \\
& \vdots & \\
\dot{r}_{n-m-1} & = & \sigma_{2} \, r_{n-m-1} 
\end{array}
\end{equation*}
where $m$ is the number of pairs of complex conjugate Floquet exponents. This corresponds to the phase-amplitude dynamics considered for instance in \cite{Llave_param,Guillamon}. The phase $\theta_{1}$ captures the periodic dynamics along the limit cycle, while other phases $\theta_j$ and amplitudes $r_j$ (with $j>1$) capture the transient dynamics toward the limit cycle. Note that the additional phases $\theta_j$ ($j>1$) are related to the dynamics of trajectories swirling in the transverse direction to the limit cycle. The limit cycle is associated with amplitude coordinates $r_j=0$ for all $j$.

The eigenfunctions $\phi_{i \omega}$ and $\phi_{\Lambda_1}$, with $\Re\{\Lambda_1\} \geq \Re\{\Lambda_j\}$ for all $j$, capture the dominant asymptotic dynamics. We can consider only the related coordinates $\theta = \angle \phi_{i \omega}$ and $r=\phi_{\Lambda_1}$ ($r=|\phi_{\Lambda_1}|$ if $\Lambda_1 \notin \mathbb{R}$), assuming that the other amplitude coordinates associated with faster dynamics are zero, i.e. $r_j=0$ for all $j\neq 1$. Denoting $\sigma \triangleq \Re\{\Lambda_1\}$, we obtain the \emph{reduced phase-amplitude dynamics}
\begin{eqnarray}
\dot{\theta} & = & \omega . \label{eq:phase_apli_1} \\
\dot{r} & = & \sigma \, r . \label{eq:phase_apli_2}
\end{eqnarray}
The phase variable $\theta$ is related to the asymptotic periodic dynamics in the longitudinal direction with respect to the limit cycle, while the amplitude variable $r$ is related to the convergent dynamics in the transverse direction. Moreover, the level sets of $\theta = \angle \phi_{i \omega}$ are the so-called \emph{isochrons} of the limit cycle \cite{Mauroy_Mezic} and the level sets of $|r|=|\phi_{\Lambda_1}|$ are the \emph{isostables} of the limit cycle \cite{MMM_isostables,Nakao_phase_amplitude}. We note that, in the case of planar systems, this is an exact (non-reduced) representation of the system in the basin of attraction of the limit cycle. For higher-dimensional systems, the phase-amplitude reduction is accurate provided that the rate of convergence of the amplitude coordinates $r_j$ (for $j>1$) is strong enough with respect to the convergence of $r_1$.

\subsection{Phase-amplitude response}

We now consider the forced limit-cycle dynamics
\begin{equation}
\label{forced_sys}
\dot{\ve{x}}  = \ve{F}(\ve{x}) + \ve{G}(\ve{x},t) \qquad \ve{x} \in \mathbb{R}^n
\end{equation}
where $\ve{G}:\mathbb{R}^n \times \mathbb{R}^+ \to \mathbb{R}^n$ is the input function. The reduced phase-amplitude dynamics are given by
\begin{eqnarray*}
\dot{\theta} & = & \nabla_x \theta \cdot (\ve{F}(\ve{x}) + \ve{G}(\ve{x},t)) = \omega + \nabla_x \theta \cdot  \ve{G}(\ve{x},t) \\
\dot{r} & = & \nabla_x r \cdot (\ve{F}(\ve{x}) + \ve{G}(\ve{x},t)) = \sigma \, r + \nabla_x r \cdot \ve{G}(\ve{x},t)
\end{eqnarray*}
where $\cdot$ denotes the inner product and $\nabla_x$ denotes the gradient with respect to the state $\ve{x}$. The phase-amplitude dynamics are given in \cite{Nakao_phase_amplitude}, where $\ve{G}$ is interpreted as a perturbation of the vector field $\ve{F}$. The gradient of $\theta$ is the \emph{phase response function} (PRF) $\ve{Z}_\theta(\theta,r)$ and the gradient of $r$ is the \emph{isostable response function} $\ve{Z}_r(\theta,r)$ (IRF). See also the definition of amplitude and phase response functions in \cite{Guillamon2}. These two functions can be expressed in terms of eigenfunctions of the Koopman operator:
\begin{equation}
\label{eq:PRS}
\begin{array}{l}
\displaystyle
\ve{Z}_\theta(\theta,r) \triangleq  \nabla_x \theta =  \nabla_x \angle \phi_{i \omega}(\ve{x}(\theta,r)) = \frac{\nabla_x \phi_{i \omega}(\ve{x}(\theta,r))}{i \, \phi_{i \omega}(\ve{x}(\theta,r))} \\
\ve{Z}_r(\theta,r) \triangleq \nabla_x r = \nabla_x \phi_{\Lambda_1}(\ve{x}(\theta,r))
\end{array}
\end{equation}
with $\ve{x}(\theta,r)$ such that $\angle \phi_{i \omega}(\ve{x}(\theta,r))= \theta$, $\phi_{\Lambda_1}(\ve{x}(\theta,r))=r$, and $\phi_{\Lambda_j}(\ve{x}(\theta,r))=0$ for all $j\neq 1$. Note that $\phi_{\Lambda_1}(\ve{x})$ should be replaced by $|\phi_{\Lambda_1}(\ve{x})|$ if $\Lambda_1 \notin \mathbb{R}$. We finally obtain the reduced dynamics
\begin{equation}
\label{eq:red_phase_ampli}
\begin{array}{l}
\dot{\theta} = \omega + \ve{Z}_\theta(\theta,r) \cdot \ve{G}(\ve{x}(\theta,r),t) . \\
\dot{r} = \sigma \, r + \ve{Z}_r(\theta,r) \cdot \ve{G}(\ve{x}(\theta,r),t) .
\end{array}
\end{equation}
For planar systems, \eqref{eq:red_phase_ampli} is not an approximation  of \eqref{forced_sys}, but an exact and equivalent representation of the dynamics.

\begin{remark}
\label{rem_coordinates}
For planar systems, the state $\ve{x}$ considered in the definition \eqref{eq:PRS} of the PRF and IRF is uniquely determined by the phase-amplitude coordinates $(\theta,r)$. For higher dimensional systems, it is clear that there is an infinity of state values $\ve{x}$ associated with the pair  $(\theta,r)$ and additional conditions $\phi_{\Lambda_j}(\ve{x})=0$ for all $j\neq 1$ must therefore be considered. However, computing all these eigenfunctions $\phi_{\Lambda_j}$ is not easy. Instead, if the system is characterized by slow-fast dynamics, one can consider the state $\ve{x}$ which lies on (an approximation of) the slow manifold of the limit cycle (where the conditions $\phi_{\Lambda_j}(\ve{x})\approx 0$, $j\neq 1$, are satisfied). 
\end{remark}

If the computation of the PRF and IRF is restricted to the limit cycle, one can further simplify the dynamics and obtain
\begin{eqnarray}
\dot{\theta} & = & \omega + \ve{Z}_\theta(\theta,0) \, \ve{G}(\ve{x}(\theta,0),t) \label{eq:phase}\\
\dot{r} & = & \sigma \, r + \ve{Z}_r(\theta,0) \, \ve{G}(\ve{x}(\theta,0),t)
\end{eqnarray}
where $\ve{Z}_\theta(\theta,0)$ is the well-known \emph{phase response curve} (PRC) \cite{Ermentrout,Izi_book} and $\ve{Z}_\theta(\theta,0)$ is the so-called \emph{isostable response curve} (IRC) defined in \cite{Wilson_isostable}. Note that \eqref{eq:phase} corresponds to classic phase reduction \cite{Brown}. Since the above phase-amplitude dynamics rely on phase and isostable response curves computed in the vicinity of the limit cycle, they are valid only locally. This is in contrast to the phase-amplitude dynamics \eqref{eq:red_phase_ampli}, which takes into account the global behavior of the system. 

\section{Computation of the eigenfunctions}
\label{sec:comput_eigen}

Phase-amplitude reduction of limit cycle dynamics is strongly connected to the spectral properties of the Koopman operator. In particular, obtaining the phase-amplitude dynamics is equivalent to computing the eigenfunctions of the Koopman operator. Our main result is to propose an efficient method to compute the dominant eigenfunctions and their gradient, and equivalently to obtain the reduced phase-amplitude dynamics \eqref{eq:red_phase_ampli}.

\subsection{Time-averaging}

An eigenfunction $\phi_{\lambda}$ can be obtained by computing the following time average of a function $f$ along the trajectories of the system \cite{Mezic,Mezic_ann_rev}:
\begin{equation}
\label{eq:time_av}
f^*_{\lambda}(\ve{x}) = \lim_{T \rightarrow \infty} \frac 1 T \int_0^T f \circ \varphi^\tau(\ve{x}) \, e^{-\lambda \tau} d\tau .
\end{equation}
If $\lambda \in i \mathbb{R}$, the time average is called \emph{Fourier average}. Otherwise, it is called \emph{Laplace average}. We observe that 
\begin{equation*}
U^t f^*_{\lambda}(\ve{x)} = \lim_{T \rightarrow \infty} \frac 1 T \int_0^T f \circ \varphi^{t+\tau}(\ve{x}) \, e^{-\lambda \tau} dt = \lim_{T \rightarrow \infty} \frac 1 T \int_0^T f \circ \varphi^{\tau}(\ve{x}) \, e^{-\lambda (\tau-t)} dt = e^{\lambda t} f^*_{\lambda}(\ve{x)}
\end{equation*}
and it follows that $f^*_{\lambda}(\ve{x})$ is an eigenfunction $\phi_{\lambda}$ (associated with the eigenvalue $\lambda$) provided that the average is not equal to zero everywhere and is a well-defined function. These conditions might be satisfied only for specific functions $f$ that we describe below.

\paragraph{Fourier averages.} In the case of a limit-cycle dynamics, the eigenfunction $\phi_{i \omega}$ can be obtained with the Fourier average $f^*_{i \omega}$ for almost all choice of $f$: the Fourier average is always well-defined and is different from zero for generic observables $f$ such that
\begin{equation*}
\int_0^{2 \pi} f \circ \varphi^t(\ve{x}) \, e^{i \omega t} \, dt  \neq 0
\end{equation*}
for $\ve{x}\in \Gamma$. We refer to \cite{Mauroy_Mezic} for more details.

\paragraph{Laplace averages.} The computation of $\phi_{\Lambda_1}$ is much more involved. According to the aforementioned results, this eigenfunction should be obtained with the Laplace average $f^*_{\Lambda_1}$. This average is non zero provided that, for all $\ve{x}\in \Gamma$, the gradient $\nabla_{x} f(\ve{x})$ is not orthogonal to the left eigenvector of the monodromy matrix associated with the Floquet exponent $\Lambda_1$ (see Section \ref{sec:num}(a)). This condition is obviously satisfied for generic functions $f$. However, selecting a function $f$ such that the average is well-defined is more delicate. Since $\Re\{\Lambda_1\}<0$, the integral \eqref{eq:time_av} converges as $T \rightarrow \infty$ only if $f \circ \varphi^\tau(\ve{x})$ tends to zero, i.e. provided that
\begin{equation}
\label{eq:cond_f}
f(\ve{x}) = 0 \quad \forall \ve{x}\in \Gamma .
\end{equation}
Since there is generally no closed-form expression of the limit cycle, finding such a function is not a trivial task.

\begin{figure}[h!]
\centering
\includegraphics[width=0.4\linewidth]{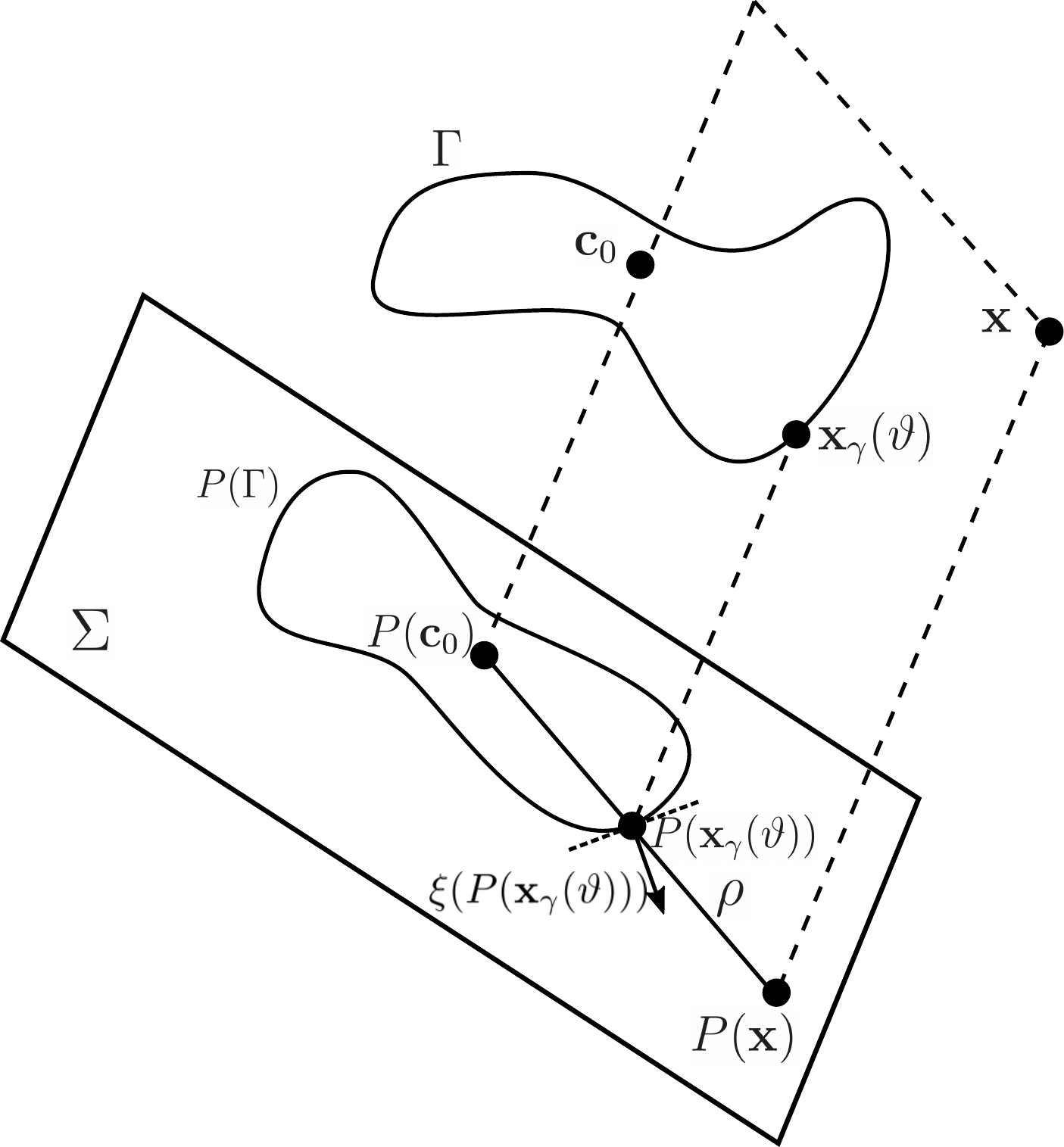}
\caption{We define polar-type coordinates $(\vartheta,\rho)$ in the plane $\Sigma$ and use them to design a function $f(\ve{x})=\rho-1$ such that $P(\ve{x})$ is assigned the coordinates $(\vartheta,\rho)$. The vector $\xi$ perpendicular to the tangent to the projected limit cycle $P(\Gamma)$ is used to compute the gradient of the eigenfunction.}
\label{fig:proj_cycle}
\end{figure}

To solve this issue, we assume that we know a parametrization of the limit cycle through its Fourier expansion
\begin{equation}
\label{eq:Fourier_coeff}
\ve{x_\gamma}:\mathbb{S} \to \Gamma, \qquad  \theta \mapsto \ve{x_\gamma}(\theta) = \sum_{k \in \mathbb{Z}} \ve{c}_k e^{i k \theta}
\end{equation}
with the phase $\theta = \omega t$. This expansion can be obtained by using a standard harmonic balance method \cite{harmonic_balance}. More details are provided in Appendix \ref{appendix}. (Note that, alternatively, one could compute the Fourier expansion of the limit cycle obtained by numerical integration of the dynamics.) Then, we suggest to choose a two-dimensional plane $\Sigma \subseteq \mathbb{R}^n$ (e.g. defined by all states equal to zero, but two) and use the following change of coordinates on that plane:
\begin{equation}
\label{eq:change_coord}
\ve{g}:\mathbb{S} \times \mathbb{R}^+ \to \Sigma , \qquad (\vartheta,\rho) \mapsto P\left( \ve{c}_0 + \rho \sum_{k \in \mathbb{Z}_*} \ve{c}_k e^{i k \vartheta}\right) =  P( \ve{c}_0) + \rho \sum_{k \in \mathbb{Z}_*} P(\ve{c}_k) e^{i k \vartheta}
\end{equation}
where $P:\mathbb{R}^n \to \Sigma$ is the orthogonal projection on $\Sigma$ (Figure \ref{fig:proj_cycle}). The map $\ve{g}$ is injective only if the interior of the curve $P(\Gamma)$ is a star set with respect to the point $P(\ve{c}_0)$. We will assume that this condition is satisfied for a well-chosen plane $\Sigma$. If it is not satisfied (e.g. in the case of possibly high-dimensional, complex limit cycles), a more sophisticated parametrization of $\Sigma$ should be considered. We leave this for future work.

Using the polar-type coordinates $(\vartheta,\rho)$ defined in \eqref{eq:change_coord}, we consider the function
\begin{equation}
\label{eq:observable_isostable}
f(\ve{x}) = \rho-1 = \frac{\|P(\ve{x})-P(\ve{c}_0)\|}{\|P(\ve{x_\gamma}(\vartheta))-P(\ve{c}_0)\|} - 1\,,
\end{equation}
where $P(\ve{x_\gamma}(\vartheta)) = \ve{g}(\vartheta,1)$, and with $\vartheta$ such that $P(\ve{x})=\ve{g}(\vartheta,\rho)$. Since $\rho=1$ for all $\ve{x}\in \Gamma$, the condition \eqref{eq:cond_f} holds and this function can be used to evaluate the Laplace average $f_{\Lambda_1}^*$.

\subsection{Gradient}

To obtain the response functions \eqref{eq:PRS} and the reduced phase-amplitude dynamics \eqref{eq:red_phase_ampli}, we have to compute the gradient of the eigenfunctions of the Koopman operator, i.e. $\nabla_x \phi_{i \omega}$ and $\nabla_x \phi_{\Lambda_1}$. This can be done through finite differences methods applied to values of the eigenfunctions, provided that these values are computed on a fine grid. Alternatively, the gradient of the eigenfunctions can also be obtained directly from time averages. Taking the gradient of \eqref{eq:time_av}, we obtain
\begin{equation*}
\nabla_x f^*_{\lambda}(\ve{x}) \cdot \delta \ve{x} = \lim_{T \rightarrow \infty} \frac 1 T \int_0^T \nabla_x(f \circ \varphi^\tau(\ve{x})) \cdot \delta \ve{x} \, e^{-\lambda \tau} d\tau = \lim_{T \rightarrow \infty} \frac 1 T \int_0^T \nabla_x f(\varphi^\tau(\ve{x})) \cdot (\ve{M}(\tau) \delta \ve{x}) \, e^{-\lambda \tau} d\tau
\end{equation*}
where $(\ve{x},\delta \ve{x}) \mapsto (\varphi^t(\ve{x}),\ve{M}(t) \delta \ve{x})$ is the flow of the prolonged system
\begin{eqnarray}
\dot{\ve{x}}  & = & \ve{F}(\ve{x}) \,, \label{eq:prol_1} \\
\dot{\delta \ve{x}} & = & \frac{\partial \ve{F}}{\partial \ve{x}} (\ve{x}) \, \delta \ve{x} \,. \label{eq:prol_2}
\end{eqnarray}
Note that $\ve{M}(t)$ is the fundamental matrix solution to $\dot{\ve{M}} = \frac{\partial \ve{F}}{\partial \ve{x}} (\ve{x}(t)) \ve{M}$, with $\ve{M}(0)=\ve{I}$. It follows that, if the time average $f_\lambda^*$ of $f$ yields an eigenfunction $\phi_{\lambda}$, then the time average $\tilde{f}^*_\lambda$ of
\begin{equation}\label{eq:f_gradient}
\tilde{f}(\ve{x,\delta x}) = \nabla_x f(\ve{x}) \cdot \ve{\delta x},
\end{equation}
computed along the trajectories of the prolonged system, yields the directional derivative of the eigenfunction $\nabla_x \phi_{\lambda} \cdot \ve{\delta x}$.
\begin{remark}
We can define the variational (or prolonged) Koopman operator $\tilde{U}^t$ associated with the prolonged system \eqref{eq:prol_1}-\eqref{eq:prol_2}, which is given by $\tilde{U}^t f(\ve{x, \delta x}) = f(\varphi^t(\ve{x}),\ve{M}(t) \delta \ve{x})$. The function $\tilde{\phi}_{\lambda}(\ve{x},\delta \ve{x}) \triangleq \nabla_x \phi_{\lambda} \cdot \ve{\delta x}$ is an eigenfunction of $\tilde{U}^t$. Indeed, we have
\begin{equation*}
\tilde{U}^t \tilde{\phi}_{\lambda}(\ve{x}, \delta x) = \tilde{U}^t \nabla_x \phi_{\lambda}(\varphi^t(\ve{x})) \cdot (\ve{M}(t) \delta \ve{x}) = \nabla_x(\phi_{\lambda} \circ \varphi^t(\ve{x})) \cdot \delta \ve{x} = e^{\lambda t} \nabla_x \phi_{\lambda}(\ve{x}) \cdot \delta \ve{x} = e^{\lambda t} \tilde{\phi}_{\lambda}(\ve{x, \delta x}) . 
\end{equation*}
\end{remark}

According to \eqref{eq:PRS}, the $j$th component of the PRF is obtained through the Fourier average
\begin{equation*}
\nabla_x f^*_{i \omega}(\ve{x}) \cdot \ve{e}_j = \lim_{T \rightarrow \infty} \frac 1 T \int_0^T \nabla_x f(\varphi^\tau(\ve{x})) \cdot (\ve{M}(\tau) \ve{e}_j) \, e^{-i \omega \tau} d\tau
\end{equation*}
for almost all $f$, where $\ve{e}_j$ is the $j$th unit vector. For the computation of the gradient of $\phi_{\Lambda_1}$, the function \eqref{eq:f_gradient} must correspond to the gradient of a function $f$ that is zero on the limit cycle. Here, we will not use the function \eqref{eq:observable_isostable} directly. Considering a function $f$ that depends on the orthogonal projection $P(\ve{x})$ (see above), we see that the gradient $\nabla_x f$ can be defined through its values on the plane $\Sigma$. Moreover, on the projected limit cycle $P(\Gamma)$, it must be perpendicular to the tangent to $P(\Gamma)$. We denote this direction by the unit vector $\xi:P(\Gamma) \to \Sigma$ (Figure \ref{fig:proj_cycle}). Using the coordinates \eqref{eq:change_coord}, we can obtain a valid candidate gradient $\nabla_x f(\ve{x})=\xi(P(\ve{x}_\gamma(\vartheta)))$ with $\vartheta$ such that $P(\ve{x})=\ve{g}(\vartheta,\rho)$ for some $\rho$. Equivalently, we define the radial projection $\pi:\Sigma \to P(\Gamma)$ such that $\pi(P(\ve{x}))=P(\ve{x_\gamma}(\vartheta))$ and, using \eqref{eq:f_gradient}, we finally obtain the function
\begin{equation*}
\tilde{f}(\ve{x},\ve{\delta x}) =\xi(P(\ve{x}_\gamma(\vartheta))) \cdot \ve{\delta x} = \xi(\pi(P(\ve{x}))) \cdot \ve{\delta x} .
\end{equation*}
According to \eqref{eq:PRS}, the $j$th component of the IRF is then obtained through the Laplace average
\begin{equation*}
\nabla_x f^*_{\Lambda_1}(\ve{x}) \cdot \ve{e}_j = \lim_{T \rightarrow \infty} \frac 1 T \int_0^T \xi(\pi(P(\varphi^\tau(\ve{x})))) \cdot (\ve{M}(\tau) \ve{e}_j) \, e^{-\Lambda_1 \tau} d\tau .
\end{equation*}

\section{Numerical computation}
\label{sec:num}

The computation of time averages (and in particular Laplace averages in the case of limit cycle dynamics) is delicate and requires some care. We provide here the following guidelines.

\setcounter{paragraph}{0}

\paragraph{Eigenvalues.}\label{sec:num_eigenval} The eigenvalues of the Koopman operator used with the time averages must be computed accurately. For the eigenvalue $i\omega$, the limit cycle frequency is obtained by the harmonic balance method (or is computed numerically by integrating the system dynamics if the harmonic balance method is not used). For the eigenvalue $\Lambda_1$, the Floquet exponent is given by $\Lambda_1 = \log(\mu_1) \, \omega/(2 \pi)$, where $\mu_1$ is the dominant eigenvalue (different from $1$) of the monodromy matrix $\ve{M}(2\pi/\omega)$, with $\ve{M}(t)$ the fundamental solution to the prolonged dynamics \eqref{eq:prol_1}-\eqref{eq:prol_2}. For planar systems, $\Lambda_1$ is also a by-product of the harmonic balance method (see Remark \ref{rem:Floquet}).

\paragraph{Fourier averages.} Fourier averages are computed with \eqref{eq:time_av}, but the integral is evaluated over a finite time horizon $\overline{T}$:
\begin{equation}
\label{eq:Fourier_av_num}
f^*_{i\omega}(\ve{x}) \approx \frac{1}{\overline{T}} \int_0^{\overline{T}} f \circ \varphi^\tau(\ve{x}) \, e^{-\lambda \tau} d\tau \qquad \overline{T} \gg 1 .
\end{equation}

\paragraph{Laplace averages.} Laplace averages are computed over a finite time horizon $\overline{T}$, which has to be finely tuned. The time horizon should be large enough for a good convergence to the limit, but also not too large so that the integrand does not blow up. In practice, we recommend to depict the value of the Laplace average (for an arbitrary initial condition) as a function of different time horizons $\overline{T}$ and to select the value $\overline{T}$ where the average reaches a constant value (before it blows up).

When the eigenvalue $\Lambda_1$ is real (e.g. case of planar systems), the Laplace average $f_{\Lambda_1}^*$ can be obtained without computing the integral in \eqref{eq:time_av}, by taking the limit
\begin{equation}
\label{eq:Laplace_av_num}
f_{\Lambda_1}^*(\ve{x}) = \lim_{T \rightarrow \infty} f(\varphi^t(\ve{x})) \, e^{-\Lambda_1 T} \approx f(\varphi^{\overline{T}}(\ve{x})) \, e^{-\Lambda_1 \overline{T}}
\end{equation}
(see \cite{MMM_isostables}). Moreover, averaging the values obtained for different values of the finite time horizon $\overline{T}_k$, $k=1\dots,K$, can provide more accurate results:
\begin{equation}
\label{eq:Laplace_av_num2}
f_{\Lambda_1}^*(\ve{x}) \approx \frac{1}{K} \sum_{k=1}^K \frac{1}{\overline{T}_k} f(\varphi^{\overline{T}_k}(\ve{x})) \, e^{-\Lambda_1 \overline{T}_k} .
\end{equation}
If $\Lambda_1$ is not real, \eqref{eq:time_av} computed over a finite time horizon yields
\begin{equation*}
f^*_{\Lambda_1}(\ve{x}) \approx \frac{1}{K} \sum_{k=1}^K \frac{1}{\overline{T}_k} \int_0^{\overline{T}_k} f \circ \varphi^\tau(\ve{x}) \, e^{-\Lambda_1 \tau} d\tau .
\end{equation*}

\paragraph{Geometric phase and bisection method.} Computing the Laplace average requires to evaluate the values of the function \eqref{eq:observable_isostable} along the trajectories, so that it is necessary to invert $\ve{x}=\ve{g}(\vartheta,\rho)$ and in particular obtain the value $\vartheta$ (the value $\vartheta$ is also needed to compute the gradient). To do so, we can use the geometric phase $\Theta(\ve{x})$ of $\ve{x}$, which we define as the signed angle \footnote{The signed angle between two vectors $\ve{v}_0$ and $\ve{v}$ can be obtained with $\textrm{atan2}(\ve{v} \times \ve{v}_0, \ve{v} \cdot \ve{v}_0)$, where $\textrm{atan2}$ is the four-quadrant inverse tangent and $\times$ denotes the vector cross-product.} between the vector $\ve{v} = \overrightarrow{P(\ve{c}_0) P(\ve{x})}$ and a reference vector $\ve{v}_0 = \overrightarrow{P(\ve{c}_0) P(\ve{x}_\gamma(0))}$. If $\ve{x}_1=\ve{g}(\vartheta,\rho_1)$ and $\ve{x}_2=\ve{g}(\vartheta,\rho_2)$, it is clear that $\Theta(\ve{x}_1)=\Theta(\ve{x}_2)$. In particular, if $\ve{x}=\ve{g}(\vartheta,\rho)$, we have the equality $\Theta(\ve{x})=\Theta(\ve{x}_\gamma(\vartheta))$ that we can use to find the value $\vartheta$ through a bisection method, exploiting the fact that $\Theta(\ve{x}_\gamma(\vartheta))$ is a monotone (increasing or decreasing) function of $\vartheta$.

\paragraph{Interpolation.} Values of Laplace and Fourier averages are computed over a uniform grid and, if needed, other values are interpolated. More details can be found in \cite{Mauroy_Mezic}. Phase and amplitude response functions are also computed over a uniform grid for specific state values. They are expressed in phase-amplitude coordinates through interpolation.

\section{Numerical examples}
\label{sec:examples}

\subsection*{Example 1: Van der Pol system}

We consider the Van der Pol system
\begin{eqnarray*}
\dot{x}_1 & = & x_2 \,. \\
\dot{x}_2 & = & x_2(1-x_1^2) - x_1 \,.
\end{eqnarray*}
The periodic orbit and the limit cycle frequency ($\omega = 0.9430$) are computed with the harmonic balance method developed in Appendix \ref{appendix}, with a Fourier series truncated to $N=40$. According to Remark \ref{rem:Floquet}, we can also obtain the Floquet exponent $\Lambda_1 = -1.059$.

\setcounter{paragraph}{0}
\paragraph{Phase-amplitude reduction.} We compute the eigenfunctions of the Koopman operator $\phi_{i \omega}$ and $\phi_{\Lambda_1}$ with the time averages evaluated on a uniform grid $100 \times 100$. The Fourier averages \eqref{eq:Fourier_av_num} are computed over a finite-time horizon $\overline{T} = 200$ with $f(\ve{x})=x_1$. The Laplace averages \eqref{eq:Laplace_av_num} are computed over a finite-time horizon $\overline{T} = 20$, with the time step specifically set to $0.1$.

\begin{figure}[h!]
\subfigure[]{\includegraphics[width=8cm]{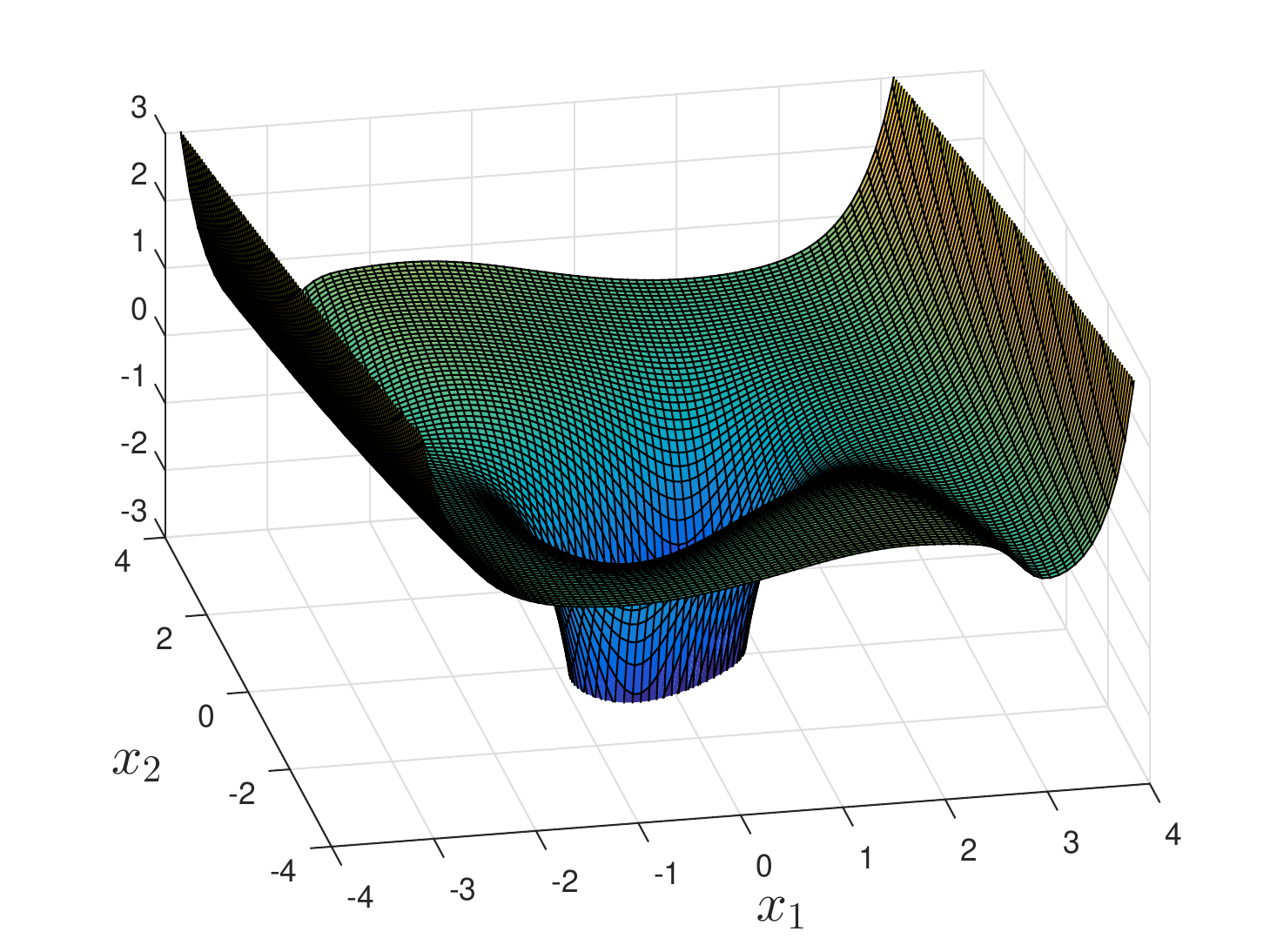}}
\subfigure[]{\includegraphics[width=8cm]{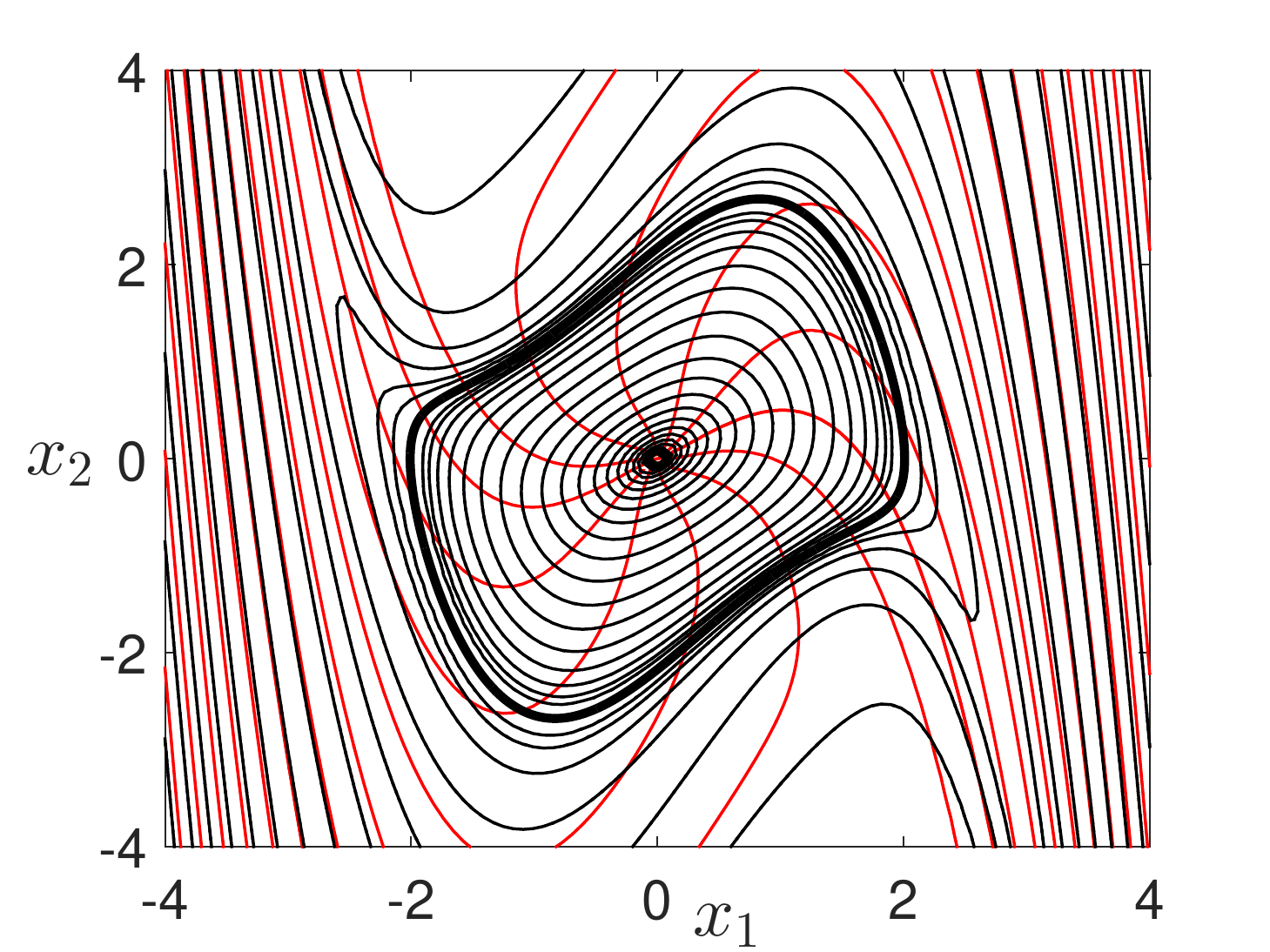}}
\caption{(a) For the Van der Pol system, the eigenfunction $\phi_{\Lambda_1}$ is computed through Laplace averages. The values associated with states on the limit cycle are in red. (b) The level sets of the Koopman operator provide exact phase amplitude-coordinates. Red: Isochrons, i.e. level sets of the eigenfunction $\phi_{i \omega}$. Black: Isostables, i.e. level sets of the Koopman eigenfunction $\phi_{\Lambda_1}$ (equally spaced on a logarithmic scale). The blue curve is the limit cycle.}
\label{fig:vdp_eigen}
\end{figure}

Figure \ref{fig:vdp_eigen}(a) shows the Koopman eigenfunction $\phi_{\Lambda_1}$. The level sets of $\phi_{i \omega}$ and $\phi_{\Lambda_1}$ (i.e. the isochrons and the isostables of the limit cycle) are depicted in \ref{fig:vdp_eigen}(b). They are the phase-amplitude coordinates of the system.

\begin{figure}[h!]
\subfigure[]{\includegraphics[width=8cm]{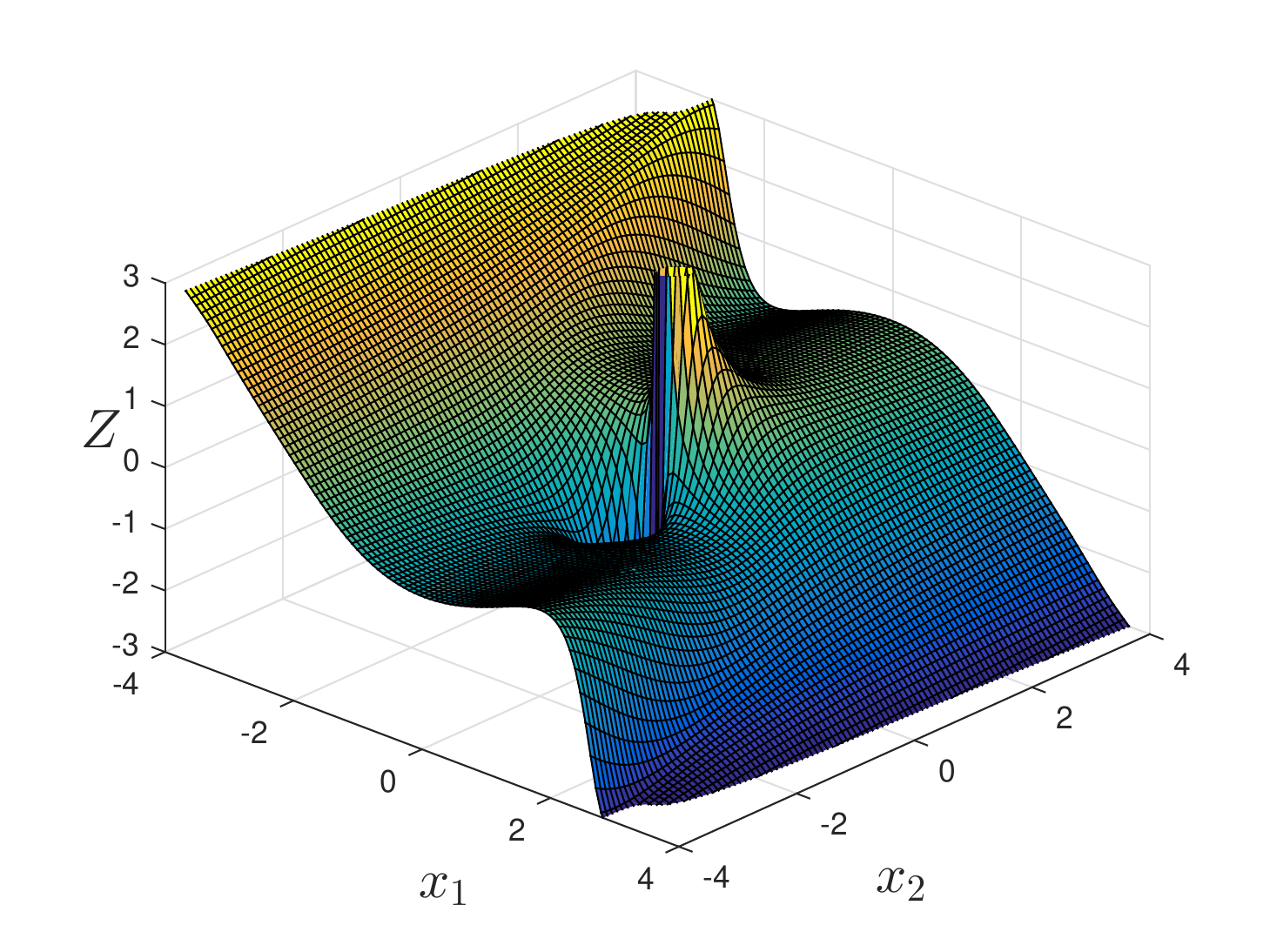}}
\subfigure[]{\includegraphics[width=8cm]{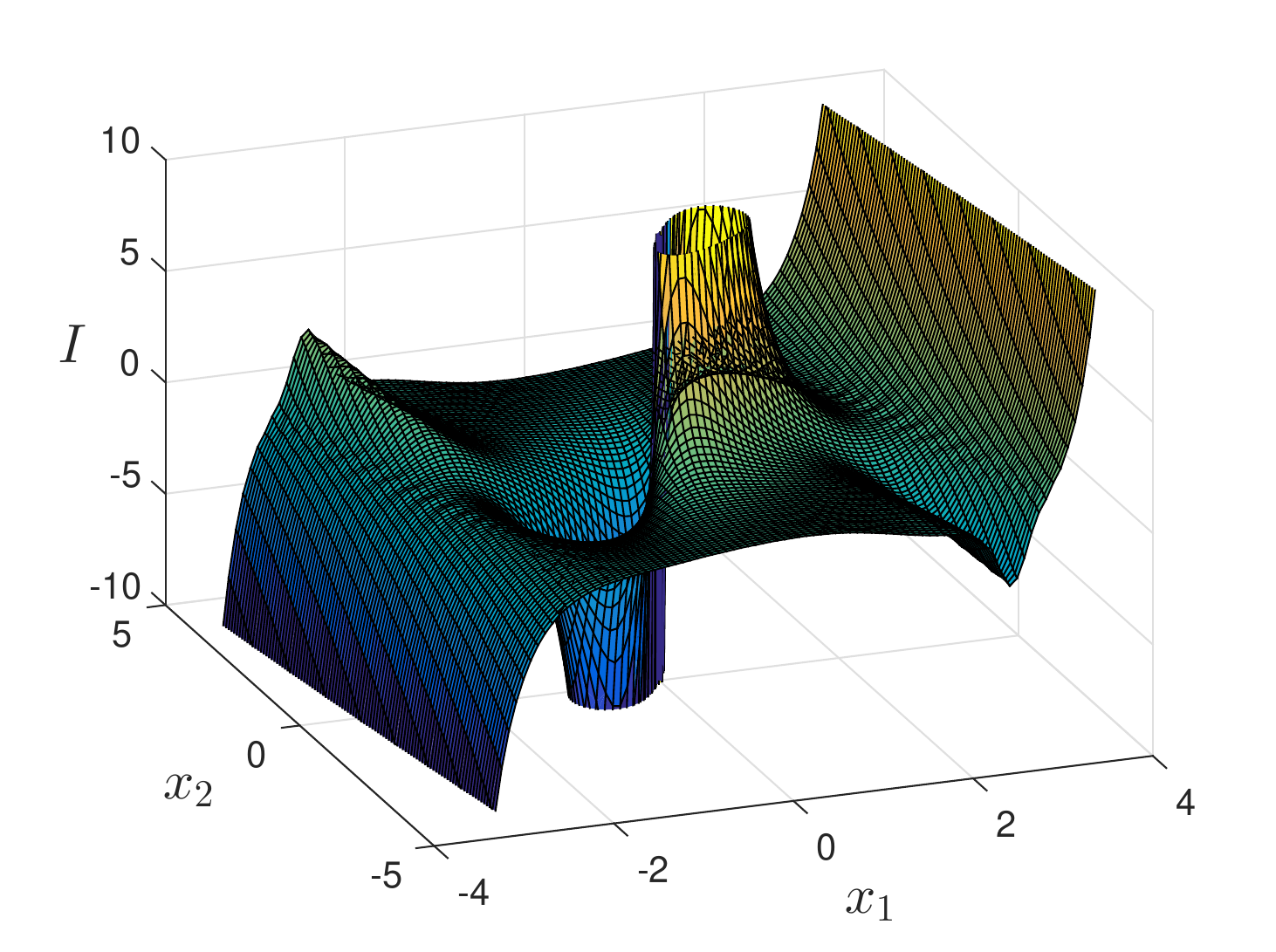}}
\caption{(a) First component of the phase response function (PRF) for the Van der Pol system. (b) First component of the isostable response function (IRF) for the Van der Pol system. For the sake of clarity, the response functions are depicted with respect to the state variables. The values associated with states on the limit cycle are in red.}
\label{fig:phase_resp}
\end{figure}

\paragraph{Phase-amplitude response.} We compute the phase and isostable response functions using Laplace averages. The first component (along $x_1$) is shown in Figure \ref{fig:phase_resp} (computations are performed with the same parameters as in the previous section).

\begin{figure}[h]
\subfigure[]{\includegraphics[width=8cm]{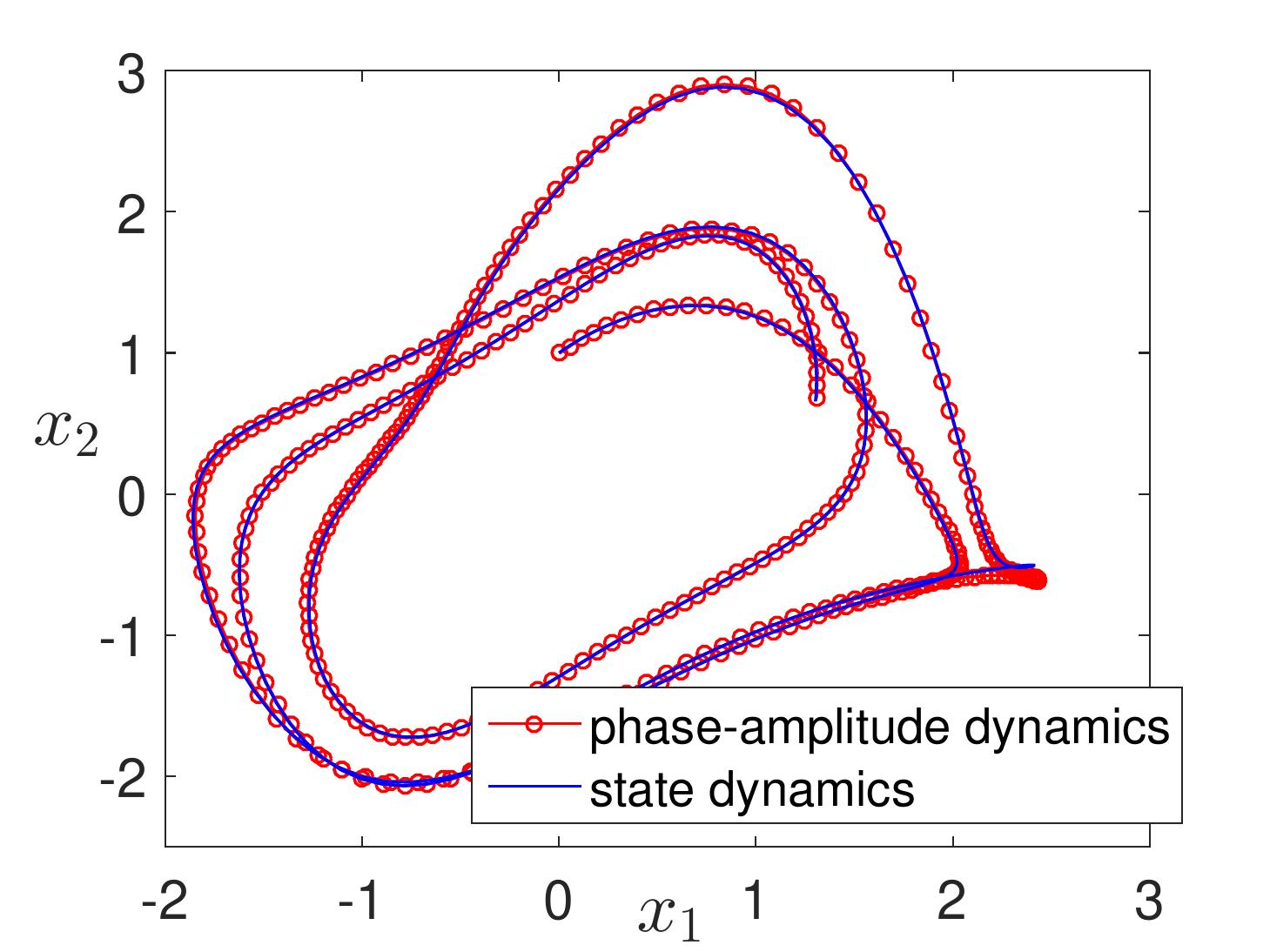}}
\subfigure[]{\includegraphics[width=8cm]{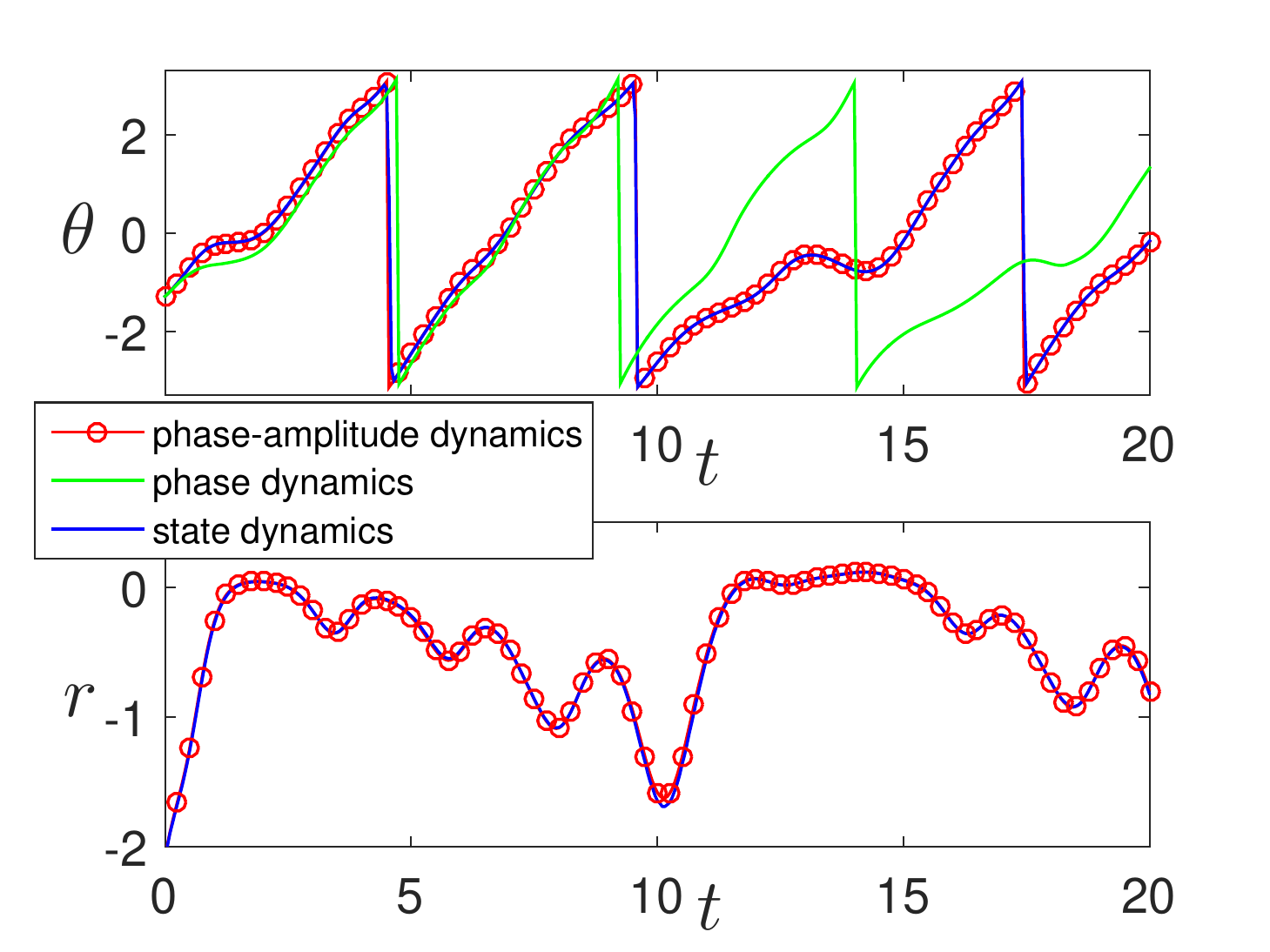}}
\caption{(a) The trajectory obtained with the phase-amplitude dynamics (red) coincides with the trajectory obtained with the original state dynamics (blue). (b) The time evolution of the phase $\theta$ and the amplitude $r$ obtained with phase-amplitude matches the evolution obtained with the state dynamics. However, classic phase reduction does not provide an accurate estimation of the evolution of the phase (green).}
\label{fig:simu_vdp}
\end{figure}

We can use the phase-amplitude dynamics \eqref{eq:red_phase_ampli} to compute the system response to an external input. In Figure \ref{fig:simu_vdp}, we compare the results obtained with the phase-amplitude dynamics and with the original state dynamics, for an input $u(t) = 0.8 \sin(1.5 t)$ applied to the first state and with the initial condition $\ve{x}=[0 , 1]$. Only very small differences are observed, which are mainly due to interpolation errors and approximations in the computation of the response functions. We also observe that classic phase reduction (obtained with the phase response curve $\ve{Z}_\theta(\theta,0)$) does not provide an accurate evolution of the phase (Figure \ref{fig:simu_vdp}(b)), since the input amplitude is not small.

\subsection*{Example 2: Three-dimensional system}

We consider a $3$-dimensional system based on the Van der Pol model:
\begin{eqnarray}
\dot{x}_1 & = & x_2 - b x_3 \label{vdp3D_1} \\
\dot{x}_2 & = & x_2(1-x_1^2) - x_1 \label{vdp3D_2} \\
\dot{x}_3 & = & a (x_1-x_3) \label{vdp3D_3}
\end{eqnarray}
with the parameters $a=2$ and $b=0.2$. The periodic orbit and its frequency ($\omega=1.1087$) are computed with the harmonic balance method, with a Fourier series truncated to $N=20$. The non zero Floquet exponents are $\Lambda_1 = -0.778$ and $\Lambda_2 = -1.843$. We note that the second Floquet exponent is related to a stronger rate of convergence, which validates the phase-amplitude reduction.

\begin{figure}[h!]
\subfigure[]{\includegraphics[width=8cm]{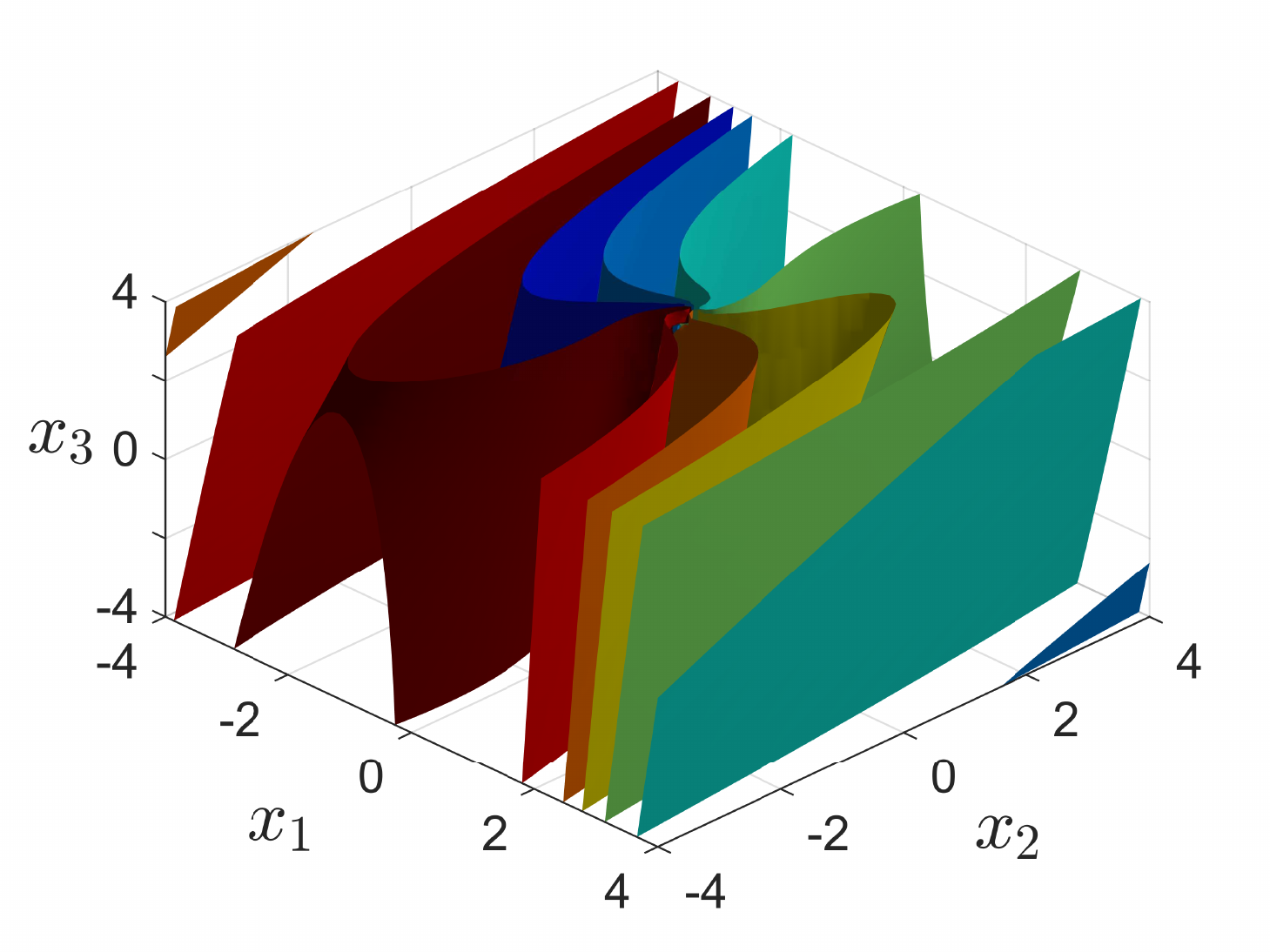}}
\subfigure[]{\includegraphics[width=8cm]{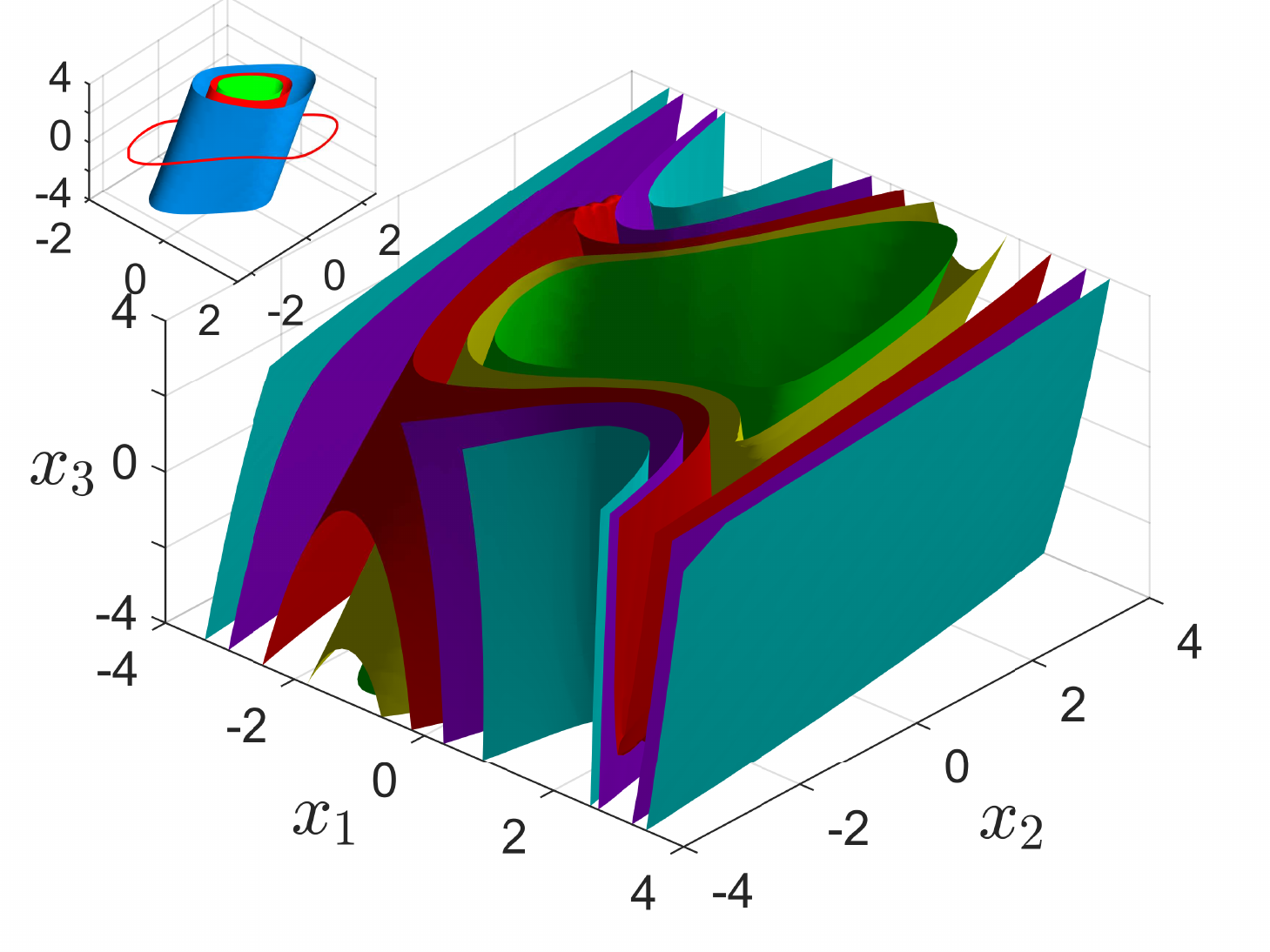}}
\caption{Isochrons and isostables of the $3$-dimensional system \eqref{vdp3D_1}-\eqref{vdp3D_3}. (a) Isochrons related to phases equally spaced by increment of $\pi/4$. (b) Isostables related to positive amplitude coordinates, equally spaced by increment of $0.1$. Inset: Isostables \guillemets{inside} the limit cycle, related to negative values equally spaced by increment of $2$. The red curve is the limit cycle.}
\label{fig:vdp_isos3D}
\end{figure}

\setcounter{paragraph}{0}
\paragraph{Phase-amplitude reduction.} We compute the eigenfunctions of the Koopman operator $\phi_{i \omega}$ and $\phi_{\Lambda_1}$ with the time averages evaluated on a uniform grid $80 \times 80 \times 80$. The Fourier averages \eqref{eq:Fourier_av_num} are computed over a finite-time horizon $\overline{T} = 200$ with $f(\ve{x})=x_1$. The Laplace averages \eqref{eq:Laplace_av_num} are computed over several finite-time horizons $\overline{T}_k \in\{24,24.1,24.2,\dots,25\}$ (with the time step specifically set to $0.1$) and by taking the average of the obtained values. We use the coordinates \eqref{eq:change_coord} where $\Sigma$ is the plane $x_3=0$. The isochrons and isostables (level sets of $\angle \phi_{i \omega}$ and $\phi_{\Lambda_1}$ respectively) are shown in Figure \ref{fig:vdp_isos3D}.

\paragraph{Phase-amplitude response.} We consider a setting that is slightly different from the case of Example 1. We suppose here that the system is forced by a train of pulses $u(t) = \epsilon \sum_n \delta(t-n\, \Delta t)$ acting on the first state. Between two pulses, the phase-amplitude dynamics is governed by \eqref{eq:phase_apli_1}-\eqref{eq:phase_apli_2}. When the system receives a pulse, the first state is instantaneously increased by $\epsilon$, which corresponds to updated phase-amplitude coordinates
\begin{equation*}
\theta \rightarrow  \theta + \Delta_\theta(\theta,r,\epsilon) \qquad r \rightarrow r + \Delta_r (\theta,r,\epsilon)
\end{equation*}
with
\begin{eqnarray*}
\Delta_\theta(\theta,r,\epsilon) & = & \angle \phi_{i \omega}(\ve{x}(\theta,r)+\epsilon \, \ve{e}_1) - \angle \phi_{i \omega}(\ve{x}(\theta,r)). \\
\Delta_r (\theta,r,\epsilon) & = & \phi_{\Lambda_1}(\ve{x}(\theta,r)+\epsilon \, \ve{e}_1) - \phi_{\Lambda_1}(\ve{x}(\theta,r)).
\end{eqnarray*}
The functions $\Delta_\theta$ and $\Delta_r$ can be seen as finite versions of the (infinitesimal) PRF and IRF, respectively, and we have 
\begin{equation*}
\lim_{\epsilon \rightarrow 0} \frac{1}{\epsilon} \Delta_\theta(\theta,r,\epsilon)) = \mathbf{Z}_\theta(\theta,r) \cdot \ve{e_1}  \qquad \lim_{\epsilon \rightarrow 0} \frac{1}{\epsilon} \Delta_r(\theta,r,\epsilon)) = \mathbf{Z}_r(\theta,r) \cdot \ve{e_1} .
\end{equation*}
In our case, the state $\ve{x}(\theta,r)$ is not fully determined by the phase-amplitude coordinates, since the system is not planar (see Remark \ref{rem_coordinates}). We consider the additional condition that the state belongs to the plane $4x-2y-5z=0$, whose approximate distance to the limit cycle is minimal (least squares minimization).

\begin{figure}[h!]
\includegraphics[width=9cm]{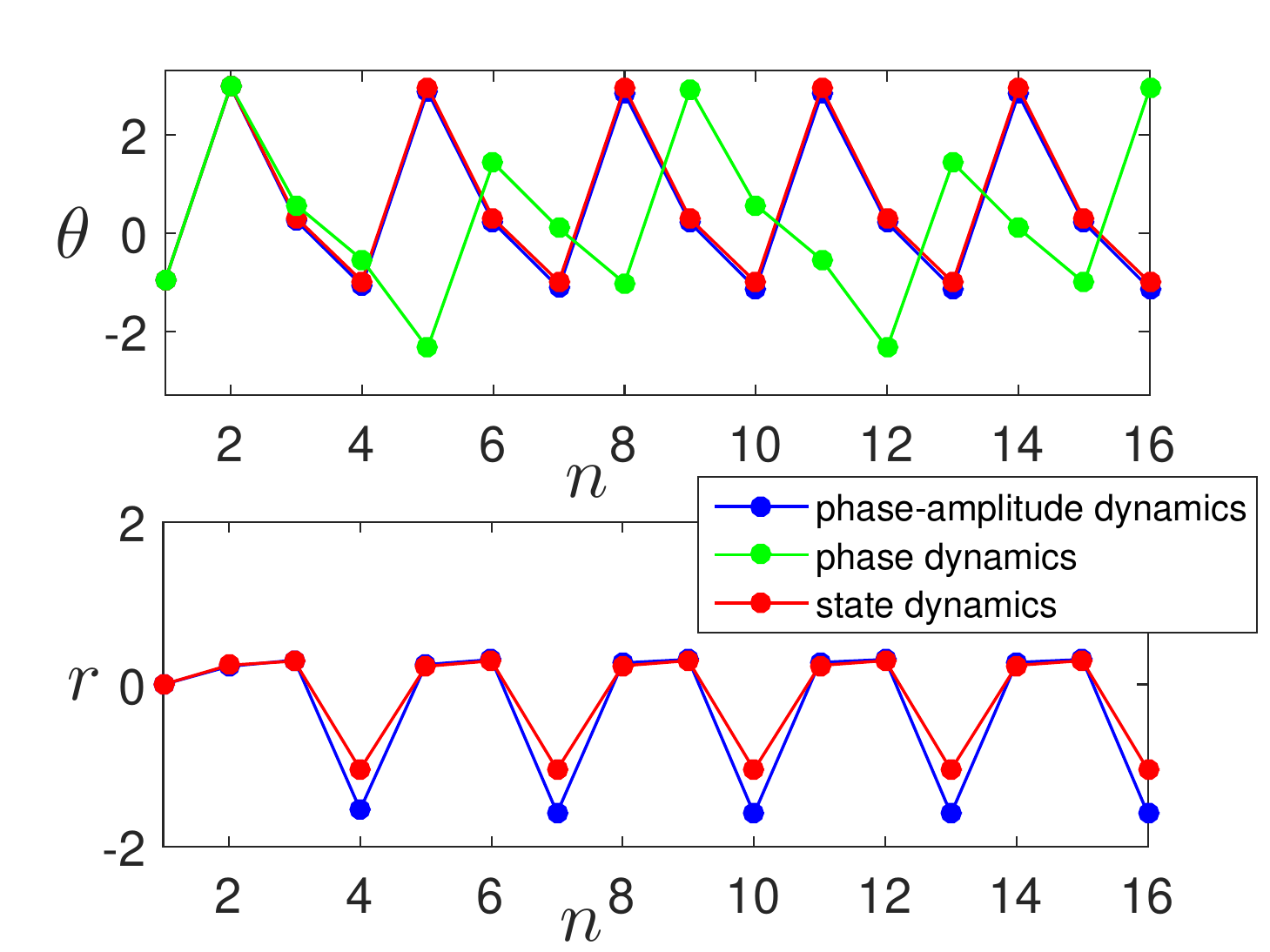}
\caption{The three-dimensional Van der Pol system is forced by a train of pulses. The phase variable (upper panel) and amplitude variables (lower panel), shown after each pulse, are correctly approximated by the phase-amplitude dynamics. However, classic phase reduction does not provide an accurate estimation of the phase dynamics (green).}
\label{fig:vdp_3D_firing}
\end{figure}

In Figure \ref{fig:vdp_3D_firing}, we compare the results obtained with the original state dynamics and with the two-dimensional phase-amplitude map
\begin{eqnarray*}
\theta[n+1] & = & \theta[n] + \omega \, \Delta t + \Delta_\theta\left( \theta[n] + \omega \, \Delta t , r[n] e^{-\Lambda_1 \, \Delta t},\epsilon\right)  \\
r[n+1] & = & r[n] e^{-\Lambda_1 \, \Delta t} + \Delta_r\left( \theta[n] + \omega \, \Delta t, r[n] e^{-\Lambda_1 \, \Delta t},\epsilon\right) 
\end{eqnarray*}
where $\theta[n]$ and $r[n]$ are the phase and amplitude values after the $n$th pulse. We consider the parameters $\epsilon=1$ and $\Delta t = 4$. Note that this map is reminiscent of the two-dimensional map considered in \cite{Guillamon2}. However, the former relies on finite phase-amplitude responses $\Delta_\theta$ and $\Delta_r$, while the latter is based on linear approximations of those responses (i.e. the infinitesimal PRF and IRF, see above). The phase-amplitude dynamics are accurate enough to provide a good approximation of the system trajectory. Main errors are observed in the amplitude coordinate and are mainly due to the fact that the phase-amplitude dynamics is a two-dimensional reduction that approximates the three-dimensional full dynamics. We verify that the one-dimensional map 
\begin{equation*}
\theta[n+1] = \theta[n] + \omega \, \Delta t + \epsilon \, \ve{Z}_\theta(\theta[n] + \omega \, \Delta t,0) \cdot \ve{e_1}
\end{equation*}
obtained from classic phase reduction does not provide an accurate evolution of the phase. We have verified that it provides accurate results when the pulse amplitude $\epsilon$ is small. We also note that, in general, the two-dimensional phase-amplitude may not be valid for periods $\Delta t$ smaller than the timescale $-1/\Lambda_2$, in which case it might be required to use a higher-dimensional map including other amplitude coordinate(s).

\subsection*{Example 3: Four-dimensional system}

We briefly show that the method is well-suited to compute the amplitude coordinates of higher-dimensional systems. We consider the $4$-dimensional Hodgkin-Huxley model (see Appendix \ref{appendix2}) which admits a limit cycle with frequency $\omega=0.429$ and a dominant Floquet exponent $\Lambda_1 = -0.178$ ($\Lambda_{2,3} =  -1.858 \pm 0.095 i$). The limit cycle is computed by numerical integration of the dynamics and its Fourier expansion is truncated to $N=150$.

\begin{figure}[h!]
\includegraphics[width=9cm]{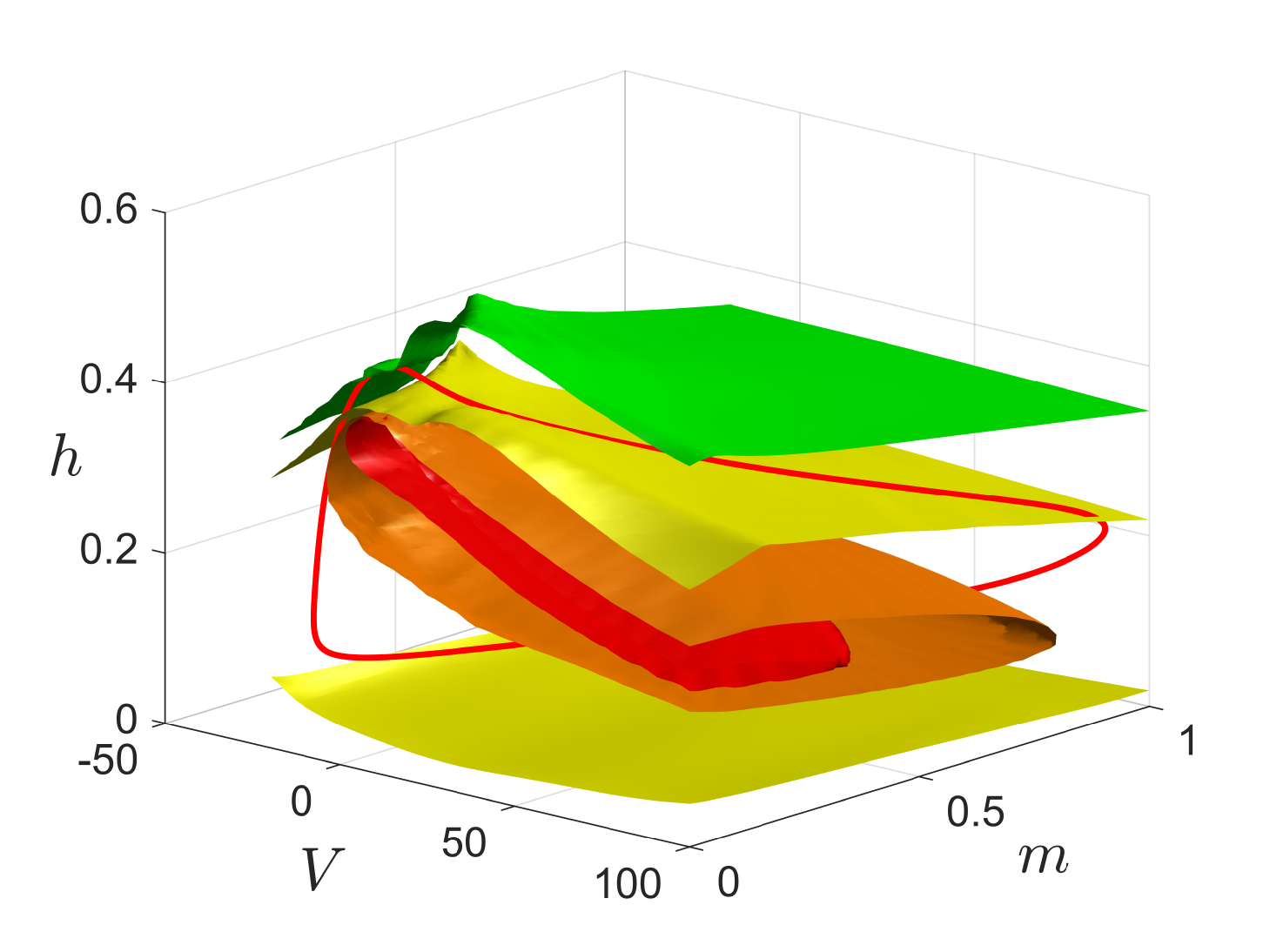}
\caption{Isostables of the $4$-dimensional Hodgkin-Huxley model. The figure shows the intersection of several isostables with the subspace $h+n=0.8$. The isostables are related to amplitude coordinates equally spaced by increment of $0.02$. The red orbit is the projection of the limit cycle on the subspace.}
\label{fig:HH}
\end{figure}

We compute the eigenfunction $\phi_{\Lambda_1}$ in the subspace $h+n=0.8$, which is known to contain the limit cycle in good approximation. The Laplace averages \eqref{eq:Laplace_av_num} are evaluated on a uniform grid $30 \times 30 \times 30$ over several finite-time horizons $\overline{T}_k \in\{70,71,72,\dots,80\}$ (with the time step specifically set to $1$) and by taking the average of the obtained values. We use the coordinates \eqref{eq:change_coord} where $\Sigma$ is the plane $(V,h)=(0,0)$. The isostables (level sets of $\phi_{\Lambda_1}$) are shown in Figure \ref{fig:HH}. Note that the computation of the isochrons (level sets of $\phi_{i \omega}$) can be found in \cite{Mauroy_Mezic}.

\section{Conclusion}
\label{sec:conclu}

In this paper, we provide an efficient method to compute phase-amplitude coordinates and responses of limit cycle dynamics, in the whole basin of attraction. In particular, we have computed the so-called isostables of limit cycle dynamics in two and three-dimensional state spaces. Our method is framed in the context of the Koopman operator, and based on the fact that phase-amplitude coordinates can be obtained by computing specific eigenfunctions of the operator. Building on previous works, we compute these eigenfunctions through Laplace averages evaluated along the system trajectories, a technique which is combined with the harmonic balance method. Our proposed method also complements previous works since it can also be used to estimate the (infinitesimal) phase and isostable response of the system. Moreover, it relies on forward integration and is therefore well-suited to compute phase and amplitude coordinates in high-dimensional spaces. 

Future work could develop a method based on a coordinate representation that is more general than the proposed polar-type coordinate representation, and therefore well-suited to (high-dimensional) systems with complex limit cycles (e.g. bursting neuron models). The method could also be extended to more general systems and used to compute the isostables and thus the response to non-local perturbations in the case of quasi-periodic tori or strange attractors, for instance. Finally, phase-amplitude coordinates could be considered to design controllers of limit cycle dynamics working with large inputs.

\section*{Acknowledgments}
I.M. work was supported in part by the ARO-MURI grant W911NF-17-1-0306 and the DARPA grant HR0011-16-C-0116.
This research used resources of the \guillemets{Plateforme Technologique de Calcul Intensif (PTCI)} located at the University of Namur, Belgium, which is supported by the F.R.S.-FNRS under the convention No.2.5020.11. The PTCI is member of the \guillemets{Consortium des Équipements de Calcul Intensif (CÉCI)}.

\appendix

\section{Approximating the limit cycle with Fourier series}
\label{appendix}

We consider the dynamical system \eqref{eq:syst} and suppose that the vector field is analytic, so that we can write
\begin{equation}
\label{eq:analytic_F}
\ve{F}(\ve{x}) = \sum_{(k_1,\dots,k_n) \in \mathbb{N}^n} \ve{F}_{k_1,\dots,k_n} \, x_1^{k_1}\cdots x_n^{k_n} .
\end{equation}
Assuming that the system admits a periodic orbit, we aim at expressing the limit cycle as a Fourier series 
\begin{equation}
\label{eq:Fourier_coeff_bis}
\ve{x_\gamma}(\theta) = \sum_{k \in \mathbb{Z}} \ve{c}_k e^{i k \theta} .
\end{equation}
To obtain the Fourier coefficients $\ve{c}_k$, we inject \eqref{eq:Fourier_coeff_bis} with $\theta = \omega t$ in the dynamics \eqref{eq:syst}, which yields
\begin{equation*}
i \omega \sum_{k \in \mathbb{Z}} k \, \ve{c}_k e^{i k \omega t}  = \ve{F} \left( \sum_{k \in \mathbb{Z}} \ve{c}_k e^{i k \omega t} \right)\,.
\end{equation*}
Using \eqref{eq:analytic_F}, we obtain
\begin{equation*}
\begin{split}
i \omega \sum_{k \in \mathbb{Z}} k \, \ve{c}_k e^{i k \omega t} = & \ve{F}_{0,\dots,0} + \sum_{\substack{(k_1,\dots,k_n) \in \mathbb{N}^n \\ k_1+\cdots+k_n>0}} \ve{F}_{k_1,\dots,k_n} \Bigg( \underbrace{\sum_{j_1 \in \mathbb{Z}} c^{(1)}_{j_1} \cdots \sum_{j_{k_1} \in \mathbb{Z}} c^{(1)}_{j_{k_1}} }_{k_1 \textrm{ times}} \Bigg) \cdots \\
& \qquad \Bigg( \underbrace{\sum_{j_{k_1+\cdots+k_{n-1}+1} \in \mathbb{Z}} c^{(n)}_{j_{k_1+\cdots+k_{n-1}+1}} \cdots \sum_{j_{k_1+\cdots+k_n} \in \mathbb{Z}} c^{(n)}_{j_{k_1+\cdots+k_n}}}_{k_n \textrm{ times}} \Bigg) \, e^{i (j_1+\cdots+j_{k_1+\cdots+k_n}) \omega t}
\end{split}
\end{equation*}
with the vectors $\ve{c}_k = (c_k^{(1)} \, \cdots \, c_k^{(n)})$. With the function
\begin{equation*}
v(k_1,\dots,k_n;l) = \begin{cases}
1 & \textrm{if } l \leq k_1 \\
r \in \{2,\dots,n\} & \textrm{if } l>k_1+\cdots+k_{r-1} \textrm{ and } l \leq k_1+\cdots+k_r
\end{cases}
\end{equation*}
we can rewrite the above expression in a more compact form:
\begin{equation*}
\begin{split}
& i \omega \sum_{k \in \mathbb{Z}} k \, \ve{c}_k e^{i k \omega t} \\
& = \ve{F}_{0,\dots,0} + \sum_{\substack{(k_1,\dots,k_n) \in \mathbb{N}^n \\ k_1+\cdots+k_n>0}} \ve{F}_{k_1,\dots,k_n} \sum_{\substack{(j_1,\dots,j_{k_1+\cdots+k_n}) \\ \in \mathbb{Z}^{k_1+\cdots+k_n}}} \prod_{l=1}^{k_1+\cdots+k_n} c_{j_l}^{v(k_1,\dots,k_n;l)} \, e^{i (j_1+\cdots+j_{k_1+\cdots+k_n}) \omega t}.
\end{split}
\end{equation*}
Finally, we obtain the set of equalities
\begin{equation}
\label{eq:equality}
i \omega k \ve{c}_k = \ve{F}_{0,\dots,0} + \sum_{\substack{(k_1,\dots,k_n) \in \mathbb{N}^n \\ k_1+\cdots+k_n>0}} \ve{F}_{k_1,\dots,k_n} \sum_{\substack{(j_1,\dots,j_{k_1+\cdots+k_n}) \in \mathbb{Z}^{k_1+\cdots+k_n} \\ j_1+\cdots+j_{k_1+\cdots+k_n}=k}} \prod_{l=1}^{k_1+\cdots+k_n} c_{j_l}^{v(k_1,\dots,k_n;l)} \,.
\end{equation}
There are as many equations as unknown Fourier coefficients (($n(2N+1)$ scalar equations if we consider the Fourier coefficients up to $k=N$). However, the frequency $\omega$ is also unknown, so that the system of equations is underdetermined. This corresponds to the fact that there are an infinity of solutions, which are all equal up to some phase lag. We can impose this phase lag, for instance by adding the constraint
\begin{equation*}
\label{eq:fix_phase}
\angle c_1^{(1)} = C
\end{equation*}
for some fixed $C\in[0,2\pi)$.

\begin{remark}[Fixed points]
A fixed point of the system yields a solution to \eqref{eq:equality} that satisfies $\ve{c}_k=0$ for $k\neq 0$. Injecting this solution in \eqref{eq:equality}, we verify that there is only one non trivial equality (for $k=0$)
\begin{equation*}
0 = \ve{F}_{0,\dots,0} + \sum_{\substack{(k_1,\dots,k_n) \in \mathbb{N}^n \\ k_1+\cdots+k_n>0}} \ve{F}_{k_1,\dots,k_n} (c_0^{(1)})^{k_1} \cdots (c_0^{(n)})^{k_n} = \ve{F}(\ve{c}_0)\,,
\end{equation*}
which holds if $\ve{c}_0$ is a fixed point of the system. This particular solution should be disregarded.
\end{remark}

\begin{remark}[Floquet exponent]
\label{rem:Floquet}
In the case of planar systems, we can express the Floquet exponent in terms of the Fourier coefficients in the expansion \eqref{eq:Fourier_coeff}. The Floquet exponent is given by
\begin{equation*}
\lambda =  \frac{1}{T} \int_0^T \textrm{div} \, \ve{F}(\ve{x_\gamma}(t)) \, dt = \frac{1}{T} \int_0^T \sum_{(k_1,k_2) \in \mathbb{N}_* \times \mathbb{N}} k_1 \, F_{k_1,k_2}^{(1)} x_1^{k_1-1} x_2^{k_2} +\sum_{(k_1,k_2) \in \mathbb{N} \times \mathbb{N}_*} k_2 \, F_{k_1,k_2}^{(2)} x_1^{k_1} x_2^{k_2-1} \, dt \,,
\end{equation*}
where $\textrm{div } \ve{F}$ denotes the divergence of $\ve{F}$ and $\ve{F}_{k_1,k_2} = (F_{k_1,k_2}^{(1)},F_{k_1,k_2}^{(2)})$. Using \eqref{eq:Fourier_coeff} and
\begin{equation*}
\frac{1}{T} \int_0^T e^{ikwt} dt =
\begin{cases}
1 & \textrm{if } k=0 \\
0 & \textrm{if } k \neq 0
\end{cases}
\end{equation*}
we obtain 
\begin{equation*}
\label{eq:Floquet}
\begin{split}
\lambda = & F_{1,0}^{(1)} + \sum_{\substack{(k_1,k_2) \in \mathbb{N}^2 \\ k_1+k_2>0}} k_1 \, F_{k_1,k_2}^{(1)}
\sum_{\substack{(j_1,\dots,j_{k_1+k_2-1}) \in \mathbb{Z}^{k_1+k_2-1} \\ j_1+\dots+j_{k_1+k_2-1}=0}} \prod_{l=1}^{k_1+k_2-1} c_{j_l}^{v(k_1-1,k_2;l)} \\
& + F_{0,1}^{(2)} + \sum_{\substack{(k_1,k_2) \in \mathbb{N}^2 \\ k_1+k_2>0}} k_2 \, F_{k_1,k_2}^{(2)}
\sum_{\substack{(j_1,\dots,j_{k_1+k_2-1}) \in \mathbb{Z}^{k_1+k_2-1} \\ j_1+\dots+j_{k_1+k_2-1}=0}} \prod_{l=1}^{k_1+k_2-1} c_{j_l}^{v(k_1,k_2-1;l)}
\end{split}.
\end{equation*}
\end{remark}

\subsection*{Numerical implementation}

We provide a few remarks and guidelines on the numerical resolution of \eqref{eq:equality}.

\begin{itemize}
\item For numerical computations, the Fourier series are truncated, i.e. we consider that $\ve{c}_k=0$ in \eqref{eq:equality} if $|k|>N$ for some chosen $N$.

\item We express \eqref{eq:equality} in terms of the real and imaginary parts of $\ve{c}_k$ and solve the obtained equations to find the (real) unknowns $\Re\{\ve{c}_k\}$ and $\Im\{\ve{c}_k\}$. Note that we can disregard the equations related to $k<0$, replacing them by the relationships
\begin{equation*}
\label{eq:conjugate}
\Re\{\ve{c}_k\} = \Re\{\ve{c}_{-k}\} \qquad \Im\{\ve{c}_k\} = -\Im\{\ve{c}_{-k}\}
\end{equation*}
since $\ve{x_\gamma}$ is real.

\item Solving \eqref{eq:equality} is a nonlinear, nonconvex problem, so that the numerical solution is likely to be inaccurate if the initial guess is too far from the true solution. To tackle this issue, we use an iterative procedure. We solve \eqref{eq:equality} for a small value $N_0$. Then the obtained solution is used as an initial guess to solve \eqref{eq:equality} with a (slightly) larger $N_1>N_0$ (for the initial guess, we assume $\ve{c}_k=0$ for $k \in  [N_0+1,N_1]$). We proceed iteratively until the error is smaller than a given threshold.

Numerical experiments suggest that, in most cases, this scheme converges and yields very accurate results. However, it might not converge when the system is high-dimensional, in particular if the limit cycle exhibits a complex geometry. In this case, another resolution scheme should be considered.
\end{itemize}

\section{Hodgkin-Huxley model}
\label{appendix2}

The Hodgkin-Huxley system is described by four states (voltage $V$ and gating variables $m$, $h$, $n$ for the (de-)activation of the ions channels Na$^+$ and K$^+$). Their dynamics are given by
\begin{eqnarray*}
\dot{V} & = & 1/C\left(-\bar{g}_{\textrm{Na}} (V-V_{\textrm{Na}}) m^3 h - \bar{g}_{\textrm{K}} (V-V_{\textrm{K}}) n^4 - g_{\textrm{L}} (V-V_L) + I_b\right) \\
\dot{m} & = & \alpha_m(V )(1 - m) - \beta_m(V )m\\
\dot{h} & = & \alpha_h(V )(1 - h) - \beta_h(V )h\\
\dot{n} & = & \alpha_n(V )(1 - n) - \beta_n(V )n
\end{eqnarray*}
with
\begin{eqnarray*}
\alpha_m(V) & = & (0.1 V - 2.5)/[1 -\exp(2.5 -0.1 V)]\\
\beta_m(V) &  = & 4 \exp(-V /18)\\
\alpha_h(V) & =  & 0.07 \exp(-V/20)\\
\beta_h(V) & =  & 1/[1 + \exp(3-0.1 V)]\\
\alpha_n(V) &  = & (0.01 V - 0.1)/[1- \exp(1-0.1 V)]\\
\beta_n(V )&  = & 0.125 \exp(-V/80)
\end{eqnarray*}
We consider the usual parameters
\begin{equation*}
\begin{split}
& V_{\textrm{Na}} = 115\, \textrm{mV},\, V_{\textrm{K}} = -12 \,\textrm{mV},\, V_{\textrm{L}} =10.6 \,\textrm{mV},\, \bar{g}_{\textrm{Na}} = 120 \,\mathrm{mS/cm^2} ,\, \bar{g}_{\mathrm{K}} = 36 \,\mathrm{mS/cm^2},\\
& \bar{g}_{\mathrm{L}} = 0.3\, \mathrm{mS/cm^2} ,\, C = 1\, \mathrm{\mu F/cm^2} , \, I_b=10 \mathrm{mA}.
\end{split}
\end{equation*}

\bibliographystyle{siam}

\end{document}